\documentclass[11pt,a4paper]{article}
\usepackage[T1]{fontenc}
\usepackage{lmodern}
\usepackage[margin=1.12in]{geometry}
\usepackage{amsmath,amssymb,amsthm,mathtools,mathrsfs,bm}
\usepackage{microtype,booktabs,enumitem,xcolor}
\usepackage[colorlinks=true,linkcolor=blue!50!black,citecolor=blue!50!black,urlcolor=blue!55!black]{hyperref}
\usepackage{cleveref}
\usepackage{url}
\allowdisplaybreaks[2]
\numberwithin{equation}{section}
\newtheorem{theorem}{Theorem}[section]
\newtheorem{proposition}[theorem]{Proposition}
\newtheorem{corollary}[theorem]{Corollary}
\newtheorem{lemma}[theorem]{Lemma}

\theoremstyle{definition}\newtheorem{definition}[theorem]{Definition}
\newtheorem{example}[theorem]{Example}
\theoremstyle{remark}\newtheorem{remark}[theorem]{Remark}
\newcommand{\Spin}{\operatorname{Spin}(7)}
\newcommand{\su}{\mathfrak{su}(2)}
\newcommand{\dd}{\mathrm d}
\newcommand{\vol}{\operatorname{vol}}
\newcommand{\pr}{\pi_7}
\newcommand{\La}{\Lambda}
\newcommand{\R}{\mathbb R}
\newcommand{\norm}[1]{\lVert #1\rVert}
\newcommand{\cof}{\operatorname{cof}}
\newcommand{\tf}{\operatorname{tf}}
\newcommand{\ea}[1]{\mathsf e_{#1}}
\newcommand{\sinf}{\sigma_\infty}
\newcommand{\DT}{\mathscr D}
\title{Finite-Energy $\Spin$-Instantons on Evolving\\ Cocalibrated $G_2$ Spaces: Damped Nahm Flows,\\ Deformation Theory and Spectral Flow\\ on the Quaternionic Heisenberg Nilmanifold}
\author{Ey\"up Yal\c{c}\i nkaya\\[2pt]
\small Scientific and Technological Research Council of T\"urkiye (T\"UB\.ITAK), Ankara, T\"urkiye\\
\small Email address: \href{mailto:eyup.yalcinkaya@tubitak.gov.tr}{\nolinkurl{eyup.yalcinkaya@tubitak.gov.tr}}}
\date{}
\begin{document}
\maketitle
\begin{abstract}
We study left-invariant $SU(2)$ $\Spin$-instantons on half-cylinders $[0,\infty)\times N^7$ whose $\Spin$ structure $\Phi=\dd t\wedge\varphi_t+\psi_t$ is induced by an evolving family of cocalibrated $G_2$ structures; for the families considered, $\Phi$ has a nonzero component in the $48$-dimensional torsion module. Explicit non-flat finite-energy examples are given on a compact quaternionic Heisenberg nilmanifold, on $S^3\times T^4$ and on $T^7$. On the Heisenberg nilmanifold, the vertical sector of the full $21$-coefficient invariant system is a damped Nahm system: it reduces to $W_\sigma=-\cof W$ on a finite interval $[0,\sinf)$ fixed by the warping. We characterise exactly which initial data give global finite-energy instantons, show that the largest admissible ball has radius $\sqrt3/\sinf$, and identify the borderline limits with an invariant $G_2$-instanton on the asymptotic slice. We then study the full, non-invariant deformation operator. Cocalibration is equivalent to self-adjointness of its slice operator, and with an Atiyah--Patodi--Singer type boundary condition the operator is Fredholm on weighted spaces, with index independent of the instanton. Using harmonic analysis on the nilmanifold and a Weitzenb\"ock identity in the Schr\"odinger representations, we determine all zero modes of the slice operators and obtain an unconditional index formula: the index counts the central characters of the lattice in an annulus fixed by the warping. For the standard integral lattice it equals $4104$.
\end{abstract}

\section{Introduction}\label{sec:intro}
The interaction between exceptional geometry and gauge theory has its origins in the theory of calibrated submanifolds. Harvey and Lawson placed the groups $G_2$ and $\Spin$ in a framework that relates special differential forms to distinguished cycles \cite{HarveyLawson1982}, and Bryant's existence results \cite{Bryant1987} started an extensive theory of special Riemannian structures in dimensions seven and eight. A parallel $G_2$ three-form or $\Spin$ four-form determines a Ricci-flat metric with restricted holonomy; allowing intrinsic torsion produces a much broader, and often more explicit, class of geometries. The relevant classifications go back to Fern\'andez--Gray for $G_2$ \cite{FernandezGray1982} and Fern\'andez for $\Spin$ \cite{Fernandez1986}, while stable forms and geometric flows provide further viewpoints \cite{Hitchin2001,Karigiannis2009}.

Gauge-theoretic analogues of four-dimensional anti-self-duality arise naturally on such manifolds \cite{DonaldsonThomas1998,Tian2000,Walpuski2017}. If $\Phi$ is the Cayley four-form of a $\Spin$ structure, a connection $A$ is a $\Spin$-instanton if
\begin{equation}\label{eq:introinst}
  F_A\in\Omega^2_{21}(M,\operatorname{ad}P),
  \qquad\text{equivalently}\qquad \pr(F_A)=0.
\end{equation}
Neither the definition of \eqref{eq:introinst} nor the ellipticity of its deformation complex requires $\dd\Phi=0$. If $\dd\Phi\neq0$, however, $\Spin$-instantons need not be Yang--Mills, and there is no topological lower bound for the curvature energy. We keep this distinction in mind throughout.

Connections with skew-symmetric torsion provide a complementary point of view \cite{FriedrichIvanov2002}. In dimension eight, every $\Spin$ structure admits a unique metric connection that preserves it and has totally skew-symmetric torsion \cite{Ivanov2004}. The intrinsic torsion of a $\Spin$ structure splits into a Lee-form component and a component in a $48$-dimensional irreducible module; the structure is locally conformally parallel exactly when the latter vanishes.

A one-parameter family $(\varphi_t,\psi_t)$ of $G_2$ structures on $N^7$ determines the $\Spin$ form
\begin{equation}\label{eq:bridgeIntro}
 \Phi=\dd t\wedge\varphi_t+\psi_t
\end{equation}
on $I\times N$. When the slices are cocalibrated, $\dd_N\psi_t=0$, one has $\dd\Phi=\dd t\wedge(\partial_t\psi_t-\dd_N\varphi_t)$. Setting this to zero gives the familiar torsion-free (Hitchin) evolution; here we deliberately do \emph{not} impose it. The cocalibrated condition has a gauge-theoretic meaning as well: it is exactly what makes the Chern--Simons functional $\int_N\operatorname{tr}(A\wedge \dd A+\frac23A^3)\wedge\psi_t$ well defined up to constants. In temporal gauge, $\Spin$-instantons on the cylinder are then its (time-dependent) gradient flow lines \cite{DonaldsonThomas1998}.

Homogeneous reductions of this kind are well established. Instanton evolution and Chern--Simons flows on cylinders and cones over seven-dimensional homogeneous spaces were studied by Haupt--Ivanova--Lechtenfeld--Popov \cite{HauptEtAl2011}, Harland--N\"olle \cite{HarlandNolle2012} and Haupt \cite{Haupt2016}. Clarke--Oliveira constructed explicit $\Spin$-instantons by evolution equations on asymptotically conical spaces \cite{ClarkeOliveira2021}. Invariant $G_2$-instantons have been studied, for instance, by Lotay--Oliveira \cite{LotayOliveira2018} and, on two-step nilpotent Lie groups, by Clarke--del Barco--Moreno \cite{ClarkeDelBarcoMoreno2025}. In the latter setting, invariant instantons on the quaternionic Heisenberg group are directly relevant to the present work (see \cref{rem:G2inst}). Ivanov--Petkov--Ugarte \cite{IvanovPetkovUgarte} investigate when the characteristic torsion connection of a $G_2$ or $\Spin$ structure is itself an instanton. Their connections live on the tangent bundle, whereas ours are auxiliary $SU(2)$ gauge fields on a trivial bundle, so the two problems are distinct.

For analytic comparison, Ghosh studies weighted $\Spin$-instanton deformation operators on asymptotically conical backgrounds and relates the index to tangential Dirac operators \cite{Ghosh2024}; our ends are asymptotically cylindrical and have nonzero intrinsic torsion. In previous work \cite{Yalcinkaya2025} the author studied deformations of instantons on \emph{locally conformal} $\Spin$ manifolds, where the torsion is carried entirely by the Lee form, the Lee-form terms cancel from the linearized equations, and the deformation theory reduces to Levi-Civita geometry. The structures considered here have $W_{48}\neq0$, so this reduction is not available; the cocalibrated cylindrical form of \cref{prop:cyl} replaces it. The nilpotent structure equations we use also appear in the rank-one solvable construction of U\u{g}uz \cite{Uguz2014}, which produces a conformally parallel $\Spin$ structure. We use only the nilpotent part of those equations: our eight-dimensional metric is a time-dependent warped product with nonvanishing $48$-dimensional torsion. No instanton result of the present paper is attributed to \cite{Uguz2014}. Cocalibrated $G_2$ structures on quaternionic Heisenberg groups and their characteristic connections were studied by Agricola--Ferreira--Storm \cite{AgricolaFerreiraStorm2015}, and relations between seven- and eight-dimensional torsion in reduction settings by Fowdar \cite{Fowdar2020}.

\paragraph{Contributions.} We state precisely what is new and what is classical.
\begin{enumerate}[label=(\arabic*),leftmargin=*]
\item On the quaternionic Heisenberg nilmanifold we write the full $21$-coefficient invariant instanton system, with all seven $G_2$-contraction equations displayed (\cref{thm:heis21}). Its horizontal curvature contains terms linear in the central coefficients, coming directly from the nilpotent brackets.
\item The vertical sector is a \emph{damped Nahm system}. An integrating factor and a time change reduce it to the $SU(2)$ Nahm equations on a finite interval $[0,\sinf)$ (\cref{thm:heiscofactor}). The Nahm equations, their reduction to the Euler top and their elliptic solutions are classical \cite{Nahm1980,Hitchin1983,Donaldson1984}. What is specific to our setting is the finite horizon $\sinf$, which is determined by the warping function.
\item We give a necessary and sufficient criterion for global existence with finite energy in terms of the classical Euler-top lifetimes (\cref{thm:heisenergy,cor:sharpT}). We show that the optimal ball of admissible data has radius $\sqrt3/\sinf$ (\cref{prop:sharpball}) and identify the borderline behaviour with convergence to an invariant $G_2$-instanton on the asymptotic slice (\cref{rem:G2inst}).
\item Curvature invariants classify gauge equivalence within the vertical family under arbitrary smooth gauge transformations (\cref{prop:heisgauge}). A short linear-ODE argument shows that solutions with integrable vertical part and decaying horizontal part are vertical (\cref{prop:heisrigidity}).
\item We set up the deformation theory of the full, non-invariant linearization (\cref{sec:deform}). We prove that cocalibration is equivalent to self-adjointness of the slice operator (\cref{prop:cyl}), that the APS-type problem is Fredholm with index constant along the vertical family (\cref{thm:fred}), and we compute the spectrum on the center-invariant sector (\cref{thm:spec0}). On the non-invariant sectors we find an explicit eigenvalue branch $4b-|\lambda|/b$ crossing zero (\cref{thm:crossing}). An exact Weitzenb\"ock identity excludes every other possible zero mode (\cref{lem:weitz,thm:centralzero}), yielding a fully unconditional lattice-dependent index formula (\cref{thm:index}) and the explicit index $4104$ for a standard integral lattice (\cref{cor:explicitindex}).
\end{enumerate}
The $S^3\times T^4$ and $T^7$ examples serve as comparisons. On $S^3\times T^4$ the vertical sector is the pull-back of the BPST family \cite{BPST1975} from the conformally cylindrical four-manifold $[0,\infty)\times S^3$. We do not claim that cofactor (Nahm) flows, homogeneous Chern--Simons flows or four-dimensional ASD reductions are new in themselves.

\paragraph{Scope.} The half-cylinder has a boundary at $t=0$ at which we impose no boundary condition. Every set of initial coefficients that extends globally therefore defines a solution. The resulting ``dimension counts'' (six, nine) are dimensions of spaces of invariant boundary data modulo gauge. They are not dimensions of moduli spaces of an elliptic boundary-value problem; the Fredholm index computed below belongs to a separately specified APS-type boundary problem. In \cref{sec:deform} we impose an APS-type boundary condition and study the resulting Fredholm problem for the full linearization.

\paragraph{Notation.} To avoid collisions we fix the following symbols: $\ea1,\ea2,\ea3$ denotes an orthonormal basis of $\su$; $\Upsilon$ is a fixed $G_2$-type three-form \eqref{eq:Q}; $H^c$ is the characteristic torsion; $\sinf$ is the cofactor horizon \eqref{eq:heisES}; $\sigma_*(W_0)$ is a maximal lifetime; $R$ denotes an element of $SO(3)$, and $O$ an element of $O(3)$ with $\det O=-1$.

\section{Cocalibrated \texorpdfstring{$G_2$}{G2} geometry and \texorpdfstring{$\Spin$}{Spin(7)} torsion}\label{sec:geom}
\subsection{Algebraic conventions}
Let $N^7$ be a smooth seven-manifold with a positively oriented orthonormal coframe $f^1,\dots,f^7$. We write $f^{i_1\cdots i_r}=f^{i_1}\wedge\cdots\wedge f^{i_r}$ and use the $G_2$ forms
\begin{align}
 \varphi={}&-f^{125}-f^{136}-f^{147}-f^{237}+f^{246}-f^{345}+f^{567},\label{eq:g2phi}\\
 \psi={}&f^{1234}-f^{1267}+f^{1357}-f^{1456}-f^{2356}-f^{2457}-f^{3467},
 \label{eq:g2psi}
\end{align}
so that $\psi=*_{7}\varphi$. The intrinsic torsion forms are defined by
\begin{align}\label{eq:g2classes}
 \dd_N\varphi&=\tau_0\psi+3\tau_1\wedge\varphi+*_{7}\tau_3,\\
 \dd_N\psi&=4\tau_1\wedge\psi+\tau_2\wedge\varphi,\nonumber
\end{align}
with $\tau_0$ a function, $\tau_1$ a one-form, $\tau_2\in\La^2_{14}$ and $\tau_3\in\La^3_{27}$. The cocalibrated condition $\dd_N\psi=0$ is equivalent to $\tau_1=\tau_2=0$, and then $\dd_N\varphi=\tau_0\psi+*_{7}\tau_3$.

Let $I\subseteq\R$ be an interval with coordinate $t$, and let $\varphi_t$ be a smooth family of positive three-forms on $N$. Put $\Phi=\dd t\wedge\varphi_t+\psi_t$ and $g_8=\dd t^2+g_t$, and orient $I\times N$ by $\dd t\wedge\vol_{7,t}$. Then
\begin{equation}\label{eq:dPhi}
 \dd_8\Phi=\dd_N\psi_t+\dd t\wedge(\dot\psi_t-\dd_N\varphi_t),
\end{equation}
so that $\dd_N\psi_t=0$ implies $\dd\Phi=\dd t\wedge(\dot\psi_t-\dd_N\varphi_t)$.

\subsection{A common family of structure equations}
Let $N$ be a seven-dimensional Lie group with a left-invariant coframe $\eta^1,\ldots,\eta^7$. Let $b:I\to(0,\infty)$ be smooth and put
\begin{equation}\label{eq:coframe}
 f^i=\eta^i\quad(1\leq i\leq4),\qquad
 f^{a+4}=b(t)\eta^{a+4}\quad(1\leq a\leq3),\qquad \beta=\frac{b'}b.
\end{equation}
Assume that on each slice
\begin{align}\label{eq:MC}
 \dd_Nf^1&=\cdots=\dd_Nf^4=0,\nonumber\\
 \dd_Nf^5&=k(f^{12}+f^{34})+m f^{67},\nonumber\\
 \dd_Nf^6&=k(f^{13}-f^{24})+m f^{75},\\
 \dd_Nf^7&=k(f^{14}+f^{23})+m f^{56},\nonumber
\end{align}
for functions $k,m$ of $t$. On $I\times N$ one has $\dd_8f^i=\dd_Nf^i$ for $i\leq4$ and
\begin{equation}\label{eq:fullMC}
 \dd_8f^j=\dd_N f^j+\beta\,\dd t\wedge f^j,\qquad j=5,6,7.
\end{equation}
Since the coframe comes from the fixed forms $\eta^i$ in \eqref{eq:coframe}, consistency of the full exterior differential also requires $k'=\beta k$ and $m'=-\beta m$; these identities hold in all three families below. The Jacobi identity ($\dd_N^2=0$) holds if and only if $km=0$; for instance, if $km\neq0$ then $\dd_N^2f^5=km(f^{137}-f^{247}-f^{146}-f^{236})\neq0$. We consider the three families
\begin{equation}\label{eq:kmcases}
 (k,m)=(-2b,0),\qquad (0,-c/b),\qquad(0,0),\qquad c>0.
\end{equation}
Define the four- and three-forms
\begin{align}\label{eq:Xi}
 \Xi={}&6f^{1234}+f^{1267}+f^{3467}-f^{1357}+f^{2457}+f^{1456}+f^{2356},\\
 \Upsilon={}&f^{125}+f^{136}+f^{147}+f^{237}-f^{246}+f^{345}+6f^{567}.\label{eq:Q}
\end{align}
They satisfy $\Xi=*_{7}\Upsilon$, $\langle\Xi,\psi\rangle=0$, $\langle \Upsilon,\varphi\rangle=0$, $\norm{\Upsilon}^2=42$, and $f^{1234}=\tfrac17(\psi+\Xi)$. Moreover $\Upsilon=h\cdot\varphi$, where $h=\operatorname{diag}(-\tfrac32,-\tfrac32,-\tfrac32,-\tfrac32,2,2,2)$ is symmetric and traceless and acts on forms as a derivation (so $h\cdot f^{ijk}=(h_i+h_j+h_k)f^{ijk}$); hence $\Upsilon\in\La^3_{27}$ and $\Xi\in\La^4_{27}$.

\begin{proposition}[Torsion of the slices and of the eight-manifold]\label{prop:geometry}
For each family in \eqref{eq:kmcases}, the forms \eqref{eq:g2phi}--\eqref{eq:g2psi} are cocalibrated on every slice, with
\begin{equation}\label{eq:tau}
 \tau_0=-\frac{6(2k-m)}7,\qquad
 \tau_3=-\frac{5k+m}{7}\,\Upsilon.
\end{equation}
The $\Spin$ structure $\Phi=\dd t\wedge\varphi_t+\psi_t$ satisfies
\begin{equation}\label{eq:torsiondecomp}
 \dd\Phi=\theta\wedge\Phi+W_{48},\qquad
 \theta=\frac{12\beta+12k-6m}{7}\,\dd t,\qquad
 W_{48}=\frac{-2\beta+5k+m}{7}\,\dd t\wedge\Xi,
\end{equation}
where $\theta\wedge\Phi$ is the Lee-form component and $W_{48}$ lies in the $48$-dimensional module. In particular, $\Phi$ is not locally conformally parallel wherever $-2\beta+5k+m\neq0$.
\end{proposition}
\begin{proof}
From \eqref{eq:MC} and the Leibniz rule one finds $\dd_N\psi=0$ and
\begin{equation}\label{eq:dphiExplicit}
 \dd_N\varphi=-6k f^{1234}+(k-m)(\Xi-6f^{1234})
 =(k-m)\Xi-\tfrac{6(2k-m)}{7}(\psi+\Xi).
\end{equation}
Since $\langle\Xi,\psi\rangle=0$ and $\norm\psi^2=7$, the $\La^4_1$-component is $\tau_0\psi$ with $\tau_0$ as in \eqref{eq:tau}, and $*_7\tau_3=-\frac{5k+m}{7}\Xi$. By \eqref{eq:coframe}, $\partial_tf^I=|I\cap\{5,6,7\}|\,\beta f^I$, so
\[
\dot\psi=2\beta(\psi-f^{1234})=\tfrac{12\beta}{7}\psi-\tfrac{2\beta}{7}\Xi.
\]
Hence
\[
\dot\psi-\dd_N\varphi=\tfrac{12\beta+12k-6m}{7}\,\psi+\tfrac{-2\beta+5k+m}{7}\,\Xi .
\]
Since $\theta\wedge\Phi=\theta(\partial_t)\,\dd t\wedge\psi$ for $\theta\in\R\,\dd t$, \eqref{eq:dPhi} gives \eqref{eq:torsiondecomp}. The $5$-forms $\alpha\wedge\Phi$, $\alpha\in T^*M$, span the $8$-dimensional module. For spatial $\alpha$, $\langle \dd t\wedge\Xi,\alpha\wedge\Phi\rangle=-\langle\Xi,\alpha\wedge\varphi\rangle=0$, because $\alpha\wedge\varphi\in\La^4_7\perp\La^4_{27}$. For $\alpha=\dd t$ the pairing is $\langle\Xi,\psi\rangle=0$. Hence $W_{48}$ lies in the complement of the Lee-form module.
\end{proof}

\begin{remark}[Characteristic connection versus gauge field]\label{rem:characteristic}
By Ivanov's theorem \cite{Ivanov2004} there is a unique metric $\Spin$-connection $\nabla^c$ with totally skew torsion $H^c$. With the orientation and Lee form fixed above, and the convention $\nabla^c=\nabla^g+\tfrac12H^c$, Ivanov's formula is
$H^c=-\delta_8\Phi-\frac76 *_8(\theta\wedge\Phi)
=*_8\dd\Phi-\frac76 *_8(\theta\wedge\Phi)$,
since $*_8\Phi=\Phi$. As $*_8(\dd t\wedge\psi)=\varphi$ and $*_8(\dd t\wedge\Xi)=\Upsilon$, \eqref{eq:torsiondecomp} gives
\begin{equation}\label{eq:charT}
 H^c=-\frac{\theta(\partial_t)}6\,\varphi+\frac{-2\beta+5k+m}{7}\,\Upsilon.
\end{equation}
Since $\varphi$ and $\Upsilon$ are orthogonal, its norm is:
\begin{equation}\label{eq:charNorm}
 \norm{H^c}^2=\frac7{36}\theta(\partial_t)^2+\frac67(-2\beta+5k+m)^2 .
\end{equation}
The connection $\nabla^c$ acts on $TM$. Our instantons are connections on a \emph{separate trivial principal $SU(2)$ bundle}, and we make no claim that $\nabla^c$ is an instanton; compare \cite{IvanovPetkovUgarte}.
\end{remark}

\section{An invariant \texorpdfstring{$SU(2)$}{SU(2)} instanton mechanism}\label{sec:mechanism}
Identify $\su$ with $(\R^3,\times)$ through an orthonormal basis $\ea1,\ea2,\ea3$ satisfying
\begin{equation}\label{eq:su2}
 [\ea1,\ea2]=\ea3,\qquad [\ea2,\ea3]=\ea1,\qquad[\ea3,\ea1]=\ea2 .
\end{equation}
Let $P=(I\times N)\times SU(2)$ and $F_A=\dd A+\tfrac12[A\wedge A]$. For a connection $A=\sum_\mu A_\mu f^\mu$, the bracket term contributes $A_\mu\times A_\nu$ to $F_{\mu\nu}$. Put $\mathcal K(\omega)=*_{8}(\Phi\wedge\omega)$. A direct computation in the orthonormal coframe $(f^1,\dots,f^7,\dd t)$ shows that $\mathcal K$ is symmetric with eigenvalues $3$ (multiplicity $7$) and $-1$ (multiplicity $21$), and that
\begin{equation}\label{eq:L27}
 \La^2_7=\operatorname{span}\{\dd t\wedge f^q+f_q\lrcorner\varphi:\ q=1,\dots,7\},\qquad
 \pr=\tfrac14(\mathrm{Id}+\mathcal K).
\end{equation}
Since $\norm{\dd t\wedge f^q+f_q\lrcorner\varphi}^2=4$ and these seven forms are mutually orthogonal, the instanton condition $\pr(F_A)=0$ is equivalent to the seven $\su$-valued equations
\begin{equation}\label{eq:G2contract}
 F_{tq}+\sum_{1\le i<j\le7}\varphi_{ijq}F_{ij}=0,\qquad q=1,\dots,7.
\end{equation}
With the nonzero components $\varphi_{125}=\varphi_{136}=\varphi_{147}=\varphi_{237}=\varphi_{345}=-1$, $\varphi_{246}=\varphi_{567}=1$ (extended skew-symmetrically), \eqref{eq:G2contract} reads
\begin{equation}\label{eq:seven}
\begin{aligned}
 &F_{t1}-F_{25}-F_{36}-F_{47}=0, &&F_{t2}+F_{15}-F_{37}+F_{46}=0,\\
 &F_{t3}+F_{16}+F_{27}-F_{45}=0, &&F_{t4}+F_{17}-F_{26}+F_{35}=0,\\
 &F_{t5}-F_{12}-F_{34}+F_{67}=0, &&F_{t6}-F_{13}+F_{24}-F_{57}=0,\\
 &F_{t7}-F_{14}-F_{23}+F_{56}=0. &&
\end{aligned}
\end{equation}

\begin{lemma}[Temporal gauge]\label{lem:temporal}
Let $A=A_t(t)\,\dd t+A_N(t)$ be a connection whose coefficients in the left-invariant coframe depend only on $t$. Then there is a smooth $g:I\to SU(2)$, depending only on $t$ and unique up to right multiplication by a constant, such that $A^g$ has no $\dd t$-component and again has $t$-dependent left-invariant coefficients.
\end{lemma}
\begin{proof}
Solve the linear ODE $g'=-A_tg$ in $SU(2)$. The coefficients of $A^g$ in the left-invariant coframe are $\operatorname{Ad}(g^{-1})$ applied to those of $A$.
\end{proof}
All invariant connections below are taken in temporal gauge. By \cref{lem:temporal}, this involves no loss of generality.

\begin{theorem}[Diagonal invariant reduction]\label{thm:reduction}
For the families \eqref{eq:MC}--\eqref{eq:kmcases}, consider
\begin{equation}\label{eq:Auvw}
 A=u(t)\ea1f^5+v(t)\ea2f^6+w(t)\ea3f^7.
\end{equation}
Then $A$ is a $\Spin$-instanton if and only if
\begin{equation}\label{eq:uvw}
 \boxed{\;u'=p\,u-vw,\quad v'=p\,v-wu,\quad w'=p\,w-uv,\qquad p=2k-m-\beta.\;}
\end{equation}
\end{theorem}
\begin{proof}
By \eqref{eq:fullMC} and \eqref{eq:su2}, the nonzero curvature components are
\begin{align}\label{eq:curv}
 &F_{12}=F_{34}=ku\,\ea1,\quad F_{13}=-F_{24}=kv\,\ea2,\quad F_{14}=F_{23}=kw\,\ea3,\nonumber\\
 &F_{67}=(vw+mu)\ea1,\quad F_{75}=(wu+mv)\ea2,\quad F_{56}=(uv+mw)\ea3,\\
 &F_{t5}=(u'+\beta u)\ea1,\quad F_{t6}=(v'+\beta v)\ea2,\quad F_{t7}=(w'+\beta w)\ea3.\nonumber
\end{align}
Each of the first four equations in \eqref{eq:seven} involves only components $F_{t i}$ ($i\le4$) and $F_{i,a+4}$, all of which vanish, so they hold identically. The last three equations give
\[
\begin{gathered}(u'+\beta u-2ku+mu+vw)\ea1=0,\qquad (v'+\beta v-2kv+mv+wu)\ea2=0,\\ (w'+\beta w-2kw+mw+uv)\ea3=0,\end{gathered}
\]
which is \eqref{eq:uvw}.
\end{proof}

\begin{theorem}[Riccati family]\label{thm:riccati}
Assume the hypotheses of \cref{thm:reduction} with $I=[0,\infty)$, and define
\begin{equation}\label{eq:EJ}
 \mathcal M(t)=\exp\Bigl(\int_0^t p(r)\,\dd r\Bigr),\qquad
 \mathcal J(t)=\int_0^t\mathcal M(r)\,\dd r.
\end{equation}
For each $C>0$, the connection
\begin{equation}\label{eq:AC}
 \boxed{A_C=z_C(t)\sum_{a=1}^3\ea{a}f^{a+4},\qquad
 z_C(t)=\frac{\mathcal M(t)}{C+\mathcal J(t)}}
\end{equation}
is a smooth $\Spin$-instanton on $[0,\infty)\times N$.
\end{theorem}
\begin{proof}
Since $\mathcal M>0$ and $\mathcal J'=\mathcal M$, we have $z_C'=pz_C-z_C^2$, and $C+\mathcal J\ge C>0$.
\end{proof}

\begin{proposition}[Non-isotropic solutions of the diagonal system]\label{prop:nondiagonal}
Write $u=\mathcal Mx$, $v=\mathcal My$, $w=\mathcal Mz$ as functions of $\mathcal J$. Then \eqref{eq:uvw} becomes
\begin{equation}\label{eq:autonomous}
 \frac{\dd x}{\dd\mathcal J}=-yz,\qquad
 \frac{\dd y}{\dd\mathcal J}=-zx,\qquad
 \frac{\dd z}{\dd\mathcal J}=-xy,
\end{equation}
whose first integrals are $x^2-y^2$ and $x^2-z^2$. For $\lambda,a_0>0$,
\begin{equation}\label{eq:hyperbolic}
 u=v=\lambda\mathcal M\operatorname{csch}\bigl(\lambda[a_0+\mathcal J]\bigr),\qquad
 w=\lambda\mathcal M\operatorname{coth}\bigl(\lambda[a_0+\mathcal J]\bigr)
\end{equation}
is a smooth instanton of the form \eqref{eq:Auvw} with $w^2-u^2=\lambda^2\mathcal M^2>0$. It is diagonal but \emph{not isotropic} ($u\neq w$).
\end{proposition}
\begin{proof}
$\mathcal M'=p\mathcal M$ and $\mathcal J'=\mathcal M$ remove the linear term. The functions $\lambda\operatorname{csch}(\lambda(a_0+\mathcal J))$ and $\lambda\coth(\lambda(a_0+\mathcal J))$ solve \eqref{eq:autonomous}, and they have no poles because $a_0+\mathcal J>0$.
\end{proof}

\begin{remark}
\eqref{eq:autonomous} is the Euler-top reduction of Nahm's equations; see \cref{rem:nahm}. The diagonal ansatz is a deliberate reduction and does not classify invariant connections.
\end{remark}

\section{Three geometrically different examples}\label{sec:examples}
In each model we use the same nonhomothetic scale
\begin{equation}\label{eq:b}
 b(t)=2-e^{-t},\qquad \beta(t)=\frac{e^{-t}}{2-e^{-t}},\qquad t\geq0,
\end{equation}
so that $1\leq b<2$, $\beta>0$ and $\beta\to0$. The metric is
\begin{equation}\label{eq:metric}
 g_8=\dd t^2+\sum_{i=1}^4(\eta^i)^2+b(t)^2\sum_{j=5}^7(\eta^j)^2.
\end{equation}
It is anisotropic and asymptotically cylindrical. The half-cylinder has a boundary at $t=0$.

\subsection{The quaternionic Heisenberg nilmanifold}\label{subsec:H}
Let $H^7$ be the simply connected two-step nilpotent Lie group with
\begin{equation}\label{eq:Heis}
 \dd\eta^1=\cdots=\dd\eta^4=0,\quad
 \dd\eta^5=-2(\eta^{12}+\eta^{34}),\quad
 \dd\eta^6=-2(\eta^{13}-\eta^{24}),\quad
 \dd\eta^7=-2(\eta^{14}+\eta^{23}).
\end{equation}
The three two-forms on the right are self-dual on $\langle\eta^1,\ldots,\eta^4\rangle$, so $H^7$ is the quaternionic Heisenberg group. Its structure constants are rational, so by Mal'cev's criterion \cite{Malcev1949} there is a cocompact lattice $\Gamma$; put $N_H=\Gamma\backslash H^7$. These are the nilpotent brackets underlying the solvable construction of \cite[\S5]{Uguz2014}. There they are extended by a derivation to a conformally parallel solvable $\Spin$ metric, whereas we use a warped product $[0,\infty)\times N_H$ with $W_{48}\neq0$.

In the moving orthonormal frame, $k=-2b$ and $m=0$. Hence
\[
 \theta_H=-\frac{12(2b-\beta)}7\,\dd t,
 \qquad W_{48,H}=-\frac{2\beta+10b}{7}\,\dd t\wedge\Xi\neq0,
\]
and the coefficient of \cref{thm:riccati} is
\begin{equation}\label{eq:EH}
 p_H=-4b-\beta,\qquad
 \mathcal M_H(t)=\frac{\exp[-8t+4(1-e^{-t})]}{b(t)},\qquad
 z_H=\frac{\mathcal M_H}{C+\int_0^t \mathcal M_H},\quad C>0.
\end{equation}
The connection $A_H=z_H\sum_a\ea{a}f^{a+4}$ is not flat, since
\begin{equation}\label{eq:Hnotflat}
 F_{12}=-2b(t)z_H(t)\,\ea1\neq0.
\end{equation}

\subsection{The product \texorpdfstring{$S^3\times T^4$}{S3 x T4}}\label{subsec:S3}
Let $N_S=SU(2)\times T^4$, with $\eta^1,\ldots,\eta^4$ closed and
\begin{equation}\label{eq:S3MC}
 \dd\eta^5=-c\eta^{67},\qquad \dd\eta^6=-c\eta^{75},\qquad \dd\eta^7=-c\eta^{56},\qquad c>0.
\end{equation}
Then $k=0$ and $m=-c/b$. Put
\begin{equation}\label{eq:S3t}
 s(t)=\int_0^t\frac{\dd r}{b(r)}=\frac12\log(2e^t-1).
\end{equation}
For $D>0$, the function
\begin{equation}\label{eq:zS}
 z_S(t)=\frac{c}{b(t)(1+D e^{-c s(t)})}
\end{equation}
satisfies $z_S'=(c/b-\beta)z_S-z_S^2$. The connection $A_S=z_S\sum_a\ea{a}f^{a+4}$ has
\begin{equation}\label{eq:Snotflat}
 F_{67}=z_S(z_S-c/b)\ea1
 =-\frac{c^2D e^{-cs}}{b^2(1+D e^{-cs})^2}\,\ea1\neq0,
 \qquad
 W_{48,S}=\frac{-2\beta-c/b}{7}\,\dd t\wedge\Xi\neq0 .
\end{equation}
As $t\to\infty$ the connection tends to the flat connection $c\sum_a\ea a\eta^{a+4}$ (the Maurer--Cartan form of $SU(2)$), not to zero.

\subsection{The seven-torus}\label{subsec:T7}
Let $N_T=T^7$ with closed $\eta^1,\ldots,\eta^7$. Then $k=m=0$ and every slice is torsion-free, while
\begin{equation}\label{eq:WT}
 \theta_T=\frac{12\beta}{7}\,\dd t,
 \qquad W_{48,T}=-\frac{2\beta}{7}\,\dd t\wedge\Xi\neq0.
\end{equation}
For $C>0$, the function $z_T=1/\bigl(b[C+s(t)]\bigr)$ satisfies $z_T'=-\beta z_T-z_T^2$, and $A_T=z_T\sum_a\ea{a}f^{a+4}$ has $F_{67}=z_T^2\ea1\neq0$. In this example the \emph{slices} have no $G_2$ torsion at all, and all the torsion of $\Phi$ comes from the evolution.

\begin{theorem}[Three nonparallel geometries with non-flat instantons]\label{thm:three}
Each of $N_H$, $N_S$, $N_T$ carries a smooth one-parameter family of cocalibrated $G_2$ structures such that $[0,\infty)\times N$ carries a $\Spin$ structure with $W_{48}\neq0$ everywhere. On these eight-manifolds, \eqref{eq:EH}, \eqref{eq:zS} and $z_T$ define smooth non-flat $SU(2)$ $\Spin$-instantons on the trivial bundle.
\end{theorem}
\begin{proof}
Combine \cref{prop:geometry}, \cref{thm:riccati} and the curvature formulas above. All displayed coefficients are nonzero because $b,\beta,c>0$.
\end{proof}

\section{Curvature energy}\label{sec:energy}
We consider
\begin{equation}\label{eq:YM}
 \mathcal E(A)=\int_{[0,\infty)\times N}\norm{F_A}^2\,\dd\vol_{g_8}.
\end{equation}
Since $\dd\Phi\neq0$, $\Spin$-instantons need not satisfy $\dd_A^*F_A=0$, and $\mathcal E$ has no topological lower bound. We therefore speak of ``curvature energy'' rather than Yang--Mills action.

\begin{lemma}[Pointwise curvature]\label{lem:Fenergy}
For $A=z\sum_a\ea{a}f^{a+4}$ with $z'=(2k-m-\beta)z-z^2$,
\begin{equation}\label{eq:Fnorm}
 \norm{F_A}^2=6k^2z^2+3(z^2+mz)^2+3\bigl((2k-m)z-z^2\bigr)^2,
\end{equation}
and $\dd\vol_{g_8}=b^3\,\dd t\wedge\dd\vol_{\eta}$.
\end{lemma}
\begin{proof}
This follows from \eqref{eq:curv}, using $z'+\beta z=(2k-m)z-z^2$.
\end{proof}

\begin{proposition}[Finite energy]\label{prop:energy}
The three connections of \cref{thm:three} have finite curvature energy.
\end{proposition}
\begin{proof}
In the Heisenberg case $z_H=O(e^{-8t})$, so \eqref{eq:Fnorm} is $O(e^{-16t})$. On $S^3\times T^4$, $z_S^2+mz_S=O(e^{-cs})$, the temporal term equals minus this, and $s=\tfrac t2+O(1)$; hence $\norm{F}^2=O(e^{-ct})$. On $T^7$, $\norm F^2=6z_T^4=O(t^{-4})$. The volume factor $b^3$ is bounded.
\end{proof}

\begin{remark}[Asymptotic torsion]
For the Heisenberg model $\beta\to0$ and $k\to-4$, so $W_{48}\to-\frac{20}{7}\dd t\wedge\Xi$ in the limiting orthonormal frame. For $S^3\times T^4$, $W_{48}\to-\frac{c}{14}\dd t\wedge\Xi$, while on $T^7$ the torsion decays. Finite gauge curvature energy is therefore compatible with non-decaying ambient torsion.
\end{remark}

\section{The full invariant system on the quaternionic Heisenberg nilmanifold}\label{sec:heisfull}
Throughout \cref{sec:heisfull,sec:heiscof,sec:heismoduli}, $N=N_H$ and $b(t)=2-e^{-t}$. Most statements only require $b$ to be smooth, bounded between positive constants, and convergent to a positive limit $b_\infty$; for \eqref{eq:b}, $b_\infty=2$. Put
\begin{equation}\label{eq:heis_s}
 s(t)=\int_0^t\frac{\dd r}{b(r)},\qquad \mathfrak a(s)=4b(s)^2 ,
\end{equation}
and consider the general invariant connection in temporal gauge
\begin{equation}\label{eq:heis_ansatz}
 A=\sum_{i=1}^4 X_i(s)\eta^i+\sum_{a=1}^3 V_a(s)\eta^{a+4},\qquad
 X_i,V_a\in\su\cong\R^3 .
\end{equation}
Write $V=[V_1\ V_2\ V_3]\in\operatorname{Mat}_{3\times3}(\R)$, $\cof V=[V_2\times V_3\ \ V_3\times V_1\ \ V_1\times V_2]$, so that $V^T\cof V=(\det V)I_3$, and
\begin{align}
 \mathcal P(X)&=\bigl[X_1\times X_2+X_3\times X_4\quad
 X_1\times X_3-X_2\times X_4\quad
 X_1\times X_4+X_2\times X_3\bigr],\label{eq:heisPX}\\
 \mathcal H(X,V)&=\bigl(
 X_2\times V_1+X_3\times V_2+X_4\times V_3,\ \
 -X_1\times V_1+X_3\times V_3-X_4\times V_2,\nonumber\\
 &\hspace{9mm}-X_1\times V_2-X_2\times V_3+X_4\times V_1,\ \
 -X_1\times V_3+X_2\times V_2-X_3\times V_1\bigr).\label{eq:heisH}
\end{align}

\begin{theorem}[Invariant Heisenberg instanton equations]\label{thm:heis21}
The connection \eqref{eq:heis_ansatz} is a $\Spin$-instanton if and only if
\begin{equation}\label{eq:heis21}
 \boxed{\quad X_s=\mathcal H(X,V),\qquad
 V_s=-\mathfrak a(s)V-\cof V+b(s)^2\mathcal P(X).\quad}
\end{equation}
The subspace $X=0$ is invariant, and on it
\begin{equation}\label{eq:heisvert}
 \boxed{\quad V_s=-\mathfrak a(s)V-\cof V.\quad}
\end{equation}
The system is equivariant under constant gauge transformations $(X_i,V)\mapsto(RX_i,RV)$, $R\in SO(3)$.
\end{theorem}
\begin{proof}
In the orthonormal coframe \eqref{eq:coframe}, $A=\sum_iX_if^i+\sum_ab^{-1}V_af^{a+4}$. By \eqref{eq:Heis} we have $\dd\eta^{5}=-2(f^{12}+f^{34})$, and similarly for $\eta^6,\eta^7$. Hence
\begin{align*}
 F_{12}&=X_1\times X_2-2V_1,&F_{34}&=X_3\times X_4-2V_1,\\
 F_{13}&=X_1\times X_3-2V_2,&F_{24}&=X_2\times X_4+2V_2,\\
 F_{14}&=X_1\times X_4-2V_3,&F_{23}&=X_2\times X_3-2V_3,
\end{align*}
together with
\[
\begin{gathered}F_{i,a+4}=b^{-1}X_i\times V_a,\qquad F_{(a+4)(d+4)}=b^{-2}V_a\times V_d,\\ F_{ti}=\partial_tX_i=b^{-1}(X_i)_s,\qquad F_{t,a+4}=b^{-2}(V_a)_s .\end{gathered}
\]
Substituting into the first four equations of \eqref{eq:seven} and multiplying by $b$ gives:
\begin{align*}
 (X_1)_s&=X_2\times V_1+X_3\times V_2+X_4\times V_3, &
 (X_2)_s&=-X_1\times V_1+X_3\times V_3-X_4\times V_2,\\
 (X_3)_s&=-X_1\times V_2-X_2\times V_3+X_4\times V_1, &
 (X_4)_s&=-X_1\times V_3+X_2\times V_2-X_3\times V_1 .
\end{align*}
The last three equations, multiplied by $b^2$, give
\begin{align*}
 (V_1)_s&=-4b^2V_1-V_2\times V_3+b^2(X_1\times X_2+X_3\times X_4),\\
 (V_2)_s&=-4b^2V_2-V_3\times V_1+b^2(X_1\times X_3-X_2\times X_4),\\
 (V_3)_s&=-4b^2V_3-V_1\times V_2+b^2(X_1\times X_4+X_2\times X_3).
\end{align*}
These are the $21$ scalar components of \eqref{eq:heis21}. Since $\mathcal H(0,V)=0$ and $\mathcal P(0)=0$, the subspace $X=0$ is invariant. Equivariance follows from $R(a\times b)=Ra\times Rb$.
\end{proof}

\begin{remark}[Gradient structure]\label{rem:gradient}
The vertical system \eqref{eq:heisvert} is the negative gradient flow, for the Frobenius metric, of the time-dependent potential
\begin{equation}\label{eq:heispot}
 \mathcal V_s(V)=2b(s)^2\norm V_F^2+\det V,
\end{equation}
because $\nabla\det V=\cof V$ holds for every matrix. This is the invariant shadow of the Chern--Simons gradient-flow structure mentioned in \cref{sec:intro}.
\end{remark}

\section{Damped Nahm reduction and the finite-energy criterion}\label{sec:heiscof}
Define
\begin{equation}\label{eq:heisES}
 \mathfrak A(s)=\int_0^s\mathfrak a,\qquad
 E(s)=e^{-\mathfrak A(s)},\qquad
 \sigma(s)=\int_0^sE,\qquad
 \sinf=\int_0^\infty E(r)\,\dd r .
\end{equation}
Since $\mathfrak a\ge4$, we have $0<\sinf\le\tfrac14$, with strict inequality for \eqref{eq:b}. Explicitly, $\mathfrak A(s(t))=4\int_0^tb$, so
\begin{equation}\label{eq:heisEt}
 E(s(t))=\exp[-8t+4(1-e^{-t})],\qquad
 \sinf=\int_0^\infty\frac{e^{-8t+4(1-e^{-t})}}{2-e^{-t}}\,\dd t\approx 0.18657 .
\end{equation}
The function $\sigma$ maps $[0,\infty)$ diffeomorphically onto $[0,\sinf)$. Since $\mathfrak a\to\mathfrak a_\infty=4b_\infty^2=16$, l'H\^opital's rule gives
\begin{equation}\label{eq:lhopital}
 \frac{E(s)}{\sinf-\sigma(s)}\longrightarrow\mathfrak a_\infty\qquad(s\to\infty).
\end{equation}

\begin{theorem}[Reduction to Nahm's equations]\label{thm:heiscofactor}
$V(s)$ solves \eqref{eq:heisvert} if and only if
\begin{equation}\label{eq:heisW}
 \boxed{\quad V(s)=E(s)\,W(\sigma(s)),\qquad
 \frac{\dd W}{\dd\sigma}=-\cof W.\quad}
\end{equation}
Along solutions,
\begin{equation}\label{eq:heisgram}
 (W^TW)_\sigma=-2(\det W)I_3,\qquad
 \tf(W^TW)\ \text{is constant},\qquad
 \tf(V^TV)=E^2\,\tf(V_0^TV_0),
\end{equation}
where $\tf(G)=G-\tfrac13\operatorname{tr}(G)I_3$. If $W(0)=R\operatorname{diag}(w_1^0,w_2^0,w_3^0)S^T$ with $R,S\in SO(3)$, then $W(\sigma)=R\operatorname{diag}(w_1,w_2,w_3)(\sigma)\,S^T$, where
\begin{equation}\label{eq:heissigned}
 w_1'=-w_2w_3,\qquad w_2'=-w_3w_1,\qquad w_3'=-w_1w_2 .
\end{equation}
\end{theorem}
\begin{proof}
Since $\cof$ is homogeneous of degree two, $V_s=-\mathfrak aEW+E^2W_\sigma$ and $-\mathfrak aV-\cof V=-\mathfrak aEW-E^2\cof W$. The Gram identity follows from $W^T\cof W=(\det W)I_3$. Every real $3\times3$ matrix has a signed singular-value decomposition with $R,S\in SO(3)$. Moreover $\cof(RDS^T)=R(\cof D)S^T$ and $\cof D=\operatorname{diag}(d_2d_3,d_3d_1,d_1d_2)$, so uniqueness of solutions gives the last claim.
\end{proof}

\begin{remark}[Nahm's equations and the Euler top]\label{rem:nahm}
Reading the columns $W_a\in\su$ with $\times=[\cdot,\cdot]$, the equation $W_\sigma=-\cof W$ becomes $\dd W_1/\dd\sigma=-[W_2,W_3]$ together with its cyclic permutations, i.e.\ the $SU(2)$ Nahm equations \cite{Nahm1980,Hitchin1983,Donaldson1984} in reversed time. Left multiplication by $SO(3)$ is the constant gauge action, and right multiplication rotates the Nahm triple. The reduction \eqref{eq:heissigned} to the Euler top, the conservation of $w_i^2-w_j^2$, and the elliptic-function solutions are classical. Our use of them is limited to computing lifetimes relative to the geometric horizon $\sinf$.
\end{remark}

Put $y(\sigma)=w_i(\sigma)^2-(w_i^0)^2$, which is independent of $i$. Then
\begin{equation}\label{eq:heisyquad}
 y'=-2w_1w_2w_3=-2\det W,\qquad (y')^2=4\prod_{i=1}^3\bigl((w_i^0)^2+y\bigr).
\end{equation}

\begin{lemma}[Blow-up profile]\label{lem:heisblowup}
If $W_\sigma=-\cof W$ has finite maximal forward lifetime $\sigma_*$, then in the signed singular-value frame of \cref{thm:heiscofactor}, for $\sigma$ near $\sigma_*$,
\begin{equation}\label{eq:blowprofile}
 |w_i(\sigma)|\,(\sigma_*-\sigma)\to1\quad(i=1,2,3),\qquad
 \det W(\sigma)<0 .
\end{equation}
\end{lemma}
\begin{proof}
The vector field is polynomial, so $\norm W\to\infty$ as $\sigma\uparrow\sigma_*$. Since the $w_i^2-w_j^2$ are constant, all $|w_i|\to\infty$ and $|w_i|/|w_j|\to1$, and eventually no $w_i$ vanishes, so the signs are constant. As $(w_i^2)'=-2\det W$ and $w_i^2\to\infty$, the determinant must eventually be negative. Then $\frac{\dd}{\dd\sigma}|w_i|=-\det W/|w_i|=|w_jw_k|=|w_i|^2(1+o(1))$, so $\frac{\dd}{\dd\sigma}|w_i|^{-1}\to-1$. Integrating from $\sigma$ to $\sigma_*$ gives \eqref{eq:blowprofile}.
\end{proof}

\begin{proposition}[Lifetimes]\label{prop:lifetime}
Let $a_i=|w_i^0|$, let $\sigma_*(W_0)\in(0,\infty]$ be the maximal forward lifetime, and put
\[
 \mathcal I(y)=\bigl((a_1^2+y)(a_2^2+y)(a_3^2+y)\bigr)^{-1/2}.
\]
\begin{enumerate}[label=\textup{(\alph*)},leftmargin=*]
\item If $\det W_0<0$, then $\displaystyle\sigma_*(W_0)=\frac12\int_0^\infty\mathcal I(y)\,\dd y$.
\item If $\det W_0>0$ and $a_{\min}=\min_ia_i$ is attained exactly once, then
\[
 \sigma_*(W_0)=\frac12\int_{-a_{\min}^2}^{0}\mathcal I(y)\,\dd y
 +\frac12\int_{-a_{\min}^2}^{\infty}\mathcal I(y)\,\dd y .
\]
If $a_{\min}$ is attained at least twice, then $\sigma_*(W_0)=\infty$.
\item If $\operatorname{rank}W_0=2$, with $a_3=0<a_1,a_2$, then $\sigma_*(W_0)=\frac12\int_0^\infty\mathcal I(y)\,\dd y$, where the singularity at $y=0$ is integrable. If $\operatorname{rank}W_0\le1$, then $W$ is constant and $\sigma_*=\infty$.
\end{enumerate}
\end{proposition}
\begin{proof}
(a) As long as no $w_i$ vanishes, $\det W$ keeps its sign. Since $w_i^2\ge a_i^2>0$ while $y$ increases, $y$ increases monotonically to $+\infty$ and $\dd\sigma=\dd y/(2\sqrt{\prod(a_i^2+y)})$; the integral converges because $\mathcal I\sim y^{-3/2}$.
(b) Now $y$ first decreases. If the minimum is unique, $\mathcal I$ has an integrable singularity $(a_{\min}^2+y)^{-1/2}$ at $-a_{\min}^2$, so $y$ reaches $-a_{\min}^2$ in finite time. There the smallest signed singular value crosses zero transversally by \eqref{eq:heissigned}, $\det W$ changes sign, and case (a) applies from $y=-a_{\min}^2$. If the minimum is repeated, the singularity is at least $(a_{\min}^2+y)^{-1}$, which is not integrable, so the turning point is reached only as $\sigma\to\infty$ and there is no blow-up.
(c) If $a_3=0$, then $(\det W)'(0)=w_1w_2w_3'=-a_1^2a_2^2<0$, so $\det W<0$ for small $\sigma>0$ and the argument of (a) applies. If the rank is at most one, then $\cof W_0=0$.
\end{proof}

\begin{proposition}[Optimal ball]\label{prop:sharpball}
For every $W_0\neq0$,
\begin{equation}\label{eq:sharpball}
 \sigma_*(W_0)\ \ge\ \frac{\sqrt3}{\norm{W_0}_F},
\end{equation}
with equality if and only if $W_0=\rho O$ with $\rho>0$ and $O\in O(3)$, $\det O=-1$. In that case $W(\sigma)=\rho O/(1-\rho\sigma)$.
\end{proposition}
\begin{proof}
Let $n=\norm W_F$. By \eqref{eq:heisgram}, $(n^2)'=-6\det W$. By the AM--GM inequality applied to the squared singular values, $|\det W|\le(n^2/3)^{3/2}$, with equality exactly when $W^TW$ is scalar. Hence $n'\le n^2/\sqrt3$ wherever $n>0$, i.e.\ $(1/n)'\ge-1/\sqrt3$. Therefore $1/n(\sigma)\ge 1/n(0)-\sigma/\sqrt3$, and blow-up cannot occur before $\sqrt3/n(0)$. If $W_0$ is not of the stated form, then either $W_0^TW_0$ is not scalar or $\det W_0\geq0$. In both cases the inequality $n'<n^2/\sqrt3$ is strict on some interval $[0,\delta)$, so $\sigma_*>\sqrt3/n(0)$. Conversely, for $W_0=\rho O$ we have $\cof(\rho O)=-\rho^2O$, the ray is invariant, $\rho'=\rho^2$, and $\sigma_*=1/\rho=\sqrt3/\norm{W_0}_F$.
\end{proof}

\begin{definition}\label{def:DT}
$\DT:=\{V_0\in\operatorname{Mat}_{3\times3}(\R):\ \sigma_*(V_0)>\sinf\}$. By continuous dependence on initial data this set is open, and it is invariant under the action $(R,S)\cdot V_0=RV_0S^T$ of $SO(3)\times SO(3)$.
\end{definition}

\begin{theorem}[Finite-energy vertical Heisenberg instantons]\label{thm:heisenergy}
The map
\begin{equation}\label{eq:heisclassmap}
 V_0\longmapsto A_{V_0}:=\sum_{a=1}^3E(s)\,W_a(\sigma(s))\,\eta^{a+4},\qquad W(0)=V_0,
\end{equation}
is a bijection from $\DT$ onto the set of smooth, globally defined, invariant $SU(2)$ $\Spin$-instantons on $[0,\infty)\times N_H$ in the vertical sector $X=0$ with finite curvature energy. Every such connection satisfies $V(s)=O(E(s))$ and tends to the trivial flat connection. Moreover:
\begin{enumerate}[label=\textup{(\roman*)},leftmargin=*]
\item $\{V_0:\norm{V_0}_F<\sqrt3/\sinf\}\subseteq\DT$, and $\sqrt3/\sinf$ is the largest radius with this property;
\item if $\sigma_*(V_0)<\sinf$, the connection becomes singular at a finite value of $s$;
\item if $\sigma_*(V_0)=\sinf$, the connection is smooth for all $s$ but $V(s)\to16\,O$ for some $O\in O(3)$ with $\det O=-1$. Its curvature does not decay, and its energy is infinite.
\end{enumerate}
For $\operatorname{rank}V_0\ge2$ the connection is irreducible; for $\operatorname{rank}V_0=1$ it is abelian (reducible).
\end{theorem}
\begin{proof}
For $X=0$, the nonzero curvature components are the horizontal ones $F_{12}=F_{34}=-2V_1$, $F_{13}=-F_{24}=-2V_2$, $F_{14}=F_{23}=-2V_3$, the central ones $F_{(a+4)(d+4)}=b^{-2}V_a\times V_d$, and the temporal ones $F_{t,a+4}=b^{-2}(V_a)_s$. Hence
\begin{equation}\label{eq:heisFnorm}
 \norm{F_A}^2=8\norm{V}^2_F+b^{-4}\norm{\cof V}^2_F
 +\norm{4V+b^{-2}\cof V}^2_F ,
\end{equation}
and $\dd\vol_{g_8}=b^4\,\dd s\,\dd\vol_\eta$ with $1\le b<2$.

If $V_0\in\DT$, then $W$ is bounded on $[0,\sinf]$, so $V=O(E)$, $\cof V=O(E^2)$, and \eqref{eq:heisFnorm} is integrable because $E\le e^{-4s}$.

Conversely, suppose the connection is smooth for all $s\ge0$ and has finite energy. Then $\int_0^\infty\norm V_F^2\,\dd s<\infty$ by \eqref{eq:heisFnorm}. Since $\sigma$ maps $[0,\infty)$ onto $[0,\sinf)$, smoothness forces $\sigma_*\ge\sinf$, and $\sigma_*<\sinf$ would give a singularity at a finite $s$, proving (ii). Suppose $\sigma_*=\sinf$. By \cref{lem:heisblowup}, $W=R\operatorname{diag}(w_i)S^T$ with $|w_i|(\sinf-\sigma)\to1$ and eventually $\det W<0$. With \eqref{eq:lhopital} this yields $V=EW\to16\,R\operatorname{diag}(\pm1,\pm1,\pm1)S^T=16\,O$ with $\det O=-1$. This contradicts $\int\norm V^2\,\dd s<\infty$ and proves (iii) and the bijection.

(i) follows from \cref{prop:sharpball}. The data $V_0=\rho O$ with $\rho=1/\sinf$ have $\sigma_*=\sinf$, so they are not in $\DT$, which shows optimality.

If $\operatorname{rank}V\ge2$ at some time, then the values $F_{12},F_{13},F_{14}$ together with $F_{67},F_{75},F_{56}$ span $\su$. Indeed, for rank two the columns of $\cof V$ are nonzero multiples of the normal to the column plane. So the stabilizer of $A$ is central. If $\operatorname{rank}V_0=1$, then $V(s)=E(s)V_0$ takes values in a fixed line of $\su$.
\end{proof}

\begin{corollary}[Effective membership test]\label{cor:sharpT}
$V_0\in\DT$ if and only if the lifetime given by \cref{prop:lifetime} in terms of the signed singular values of $V_0$ exceeds $\sinf\approx0.18657$.
\end{corollary}

\begin{remark}[Borderline trajectories and an invariant $G_2$-instanton]\label{rem:G2inst}
On the asymptotic slice ($b\equiv2$), the nonzero equilibria of $V_s=-16V-\cof V$ are exactly the matrices $V=16\,O$ with $\det O=-1$. A time-independent solution of \eqref{eq:seven} is a connection on the slice with $F\in\La^2_{14}$, i.e.\ a $G_2$-instanton. Hence
\[
 A_\infty=16\sum_{a=1}^3(O\ea a)\,\eta^{a+4},\qquad \det O=-1,
\]
is an irreducible left-invariant $G_2$-instanton on $(N_H,\varphi|_{b=2})$. It is not flat, since $F_{12}=-32\,O\ea1$. By \cref{thm:heisenergy}, the initial data with $\sigma_*(V_0)=\sinf$ are exactly those whose trajectories are globally defined and converge to one of these instantons. We expect $A_\infty$ to be among the invariant instantons on two-step nilpotent groups classified in \cite{ClarkeDelBarcoMoreno2025}, up to normalisation, but we have not checked the identification.
\end{remark}

\begin{example}[Anisotropic irreducible data]
Let $W_0^{\pm}=\operatorname{diag}\bigl(\tfrac15,\tfrac3{10},\pm\tfrac25\bigr)$. Then $\norm{W_0^\pm}_F=\sqrt{29}/10\approx0.54<\sqrt3/\sinf\approx9.28$, so $W_0^\pm\in\DT$. Both give full-rank, irreducible, anisotropic finite-energy instantons with the same Gram matrix and opposite determinants. By \cref{prop:heisgauge} they are not gauge equivalent. Replacing $W_0^\pm$ by $W_0^\pm S^T$ with $S\in SO(3)$ gives non-diagonal coefficient matrices.
\end{example}

\begin{remark}[Explicit isotropic branches]\label{rem:heisscalar}
For $R\in SO(3)$ and $O\in O(3)$ with $\det O=-1$,
\begin{equation}\label{eq:heissignbranches}
 V_+(s)=\frac{E(s)}{C+\sigma(s)}R\quad(C>0),\qquad
 V_-(s)=\frac{E(s)}{C-\sigma(s)}O\quad(C>\sinf)
\end{equation}
are finite-energy, non-flat and mutually gauge-inequivalent. For $V_-$, the value $C=\sinf$ is the borderline case of \cref{thm:heisenergy}(iii), and $C<\sinf$ gives a singularity at finite $s$. Rank-one data have $\cof V_0=0$, so $V=EV_0$ has finite energy for every $V_0$; hence $\DT$ is unbounded and is not a ball.
\end{remark}

\section{Gauge classes and horizontal decay}\label{sec:heismoduli}
We now allow arbitrary smooth gauge transformations $g:[0,\infty)\times N_H\to SU(2)$, not necessarily invariant.

\begin{proposition}[Gauge classes in the vertical family]\label{prop:heisgauge}
Two connections $A_{V_0}$ and $A_{\widetilde V_0}$ of \cref{thm:heisenergy} are gauge equivalent if and only if
\begin{equation}\label{eq:heisorbit}
 V_0^TV_0=\widetilde V_0^T\widetilde V_0\qquad\text{and}\qquad \det V_0=\det\widetilde V_0,
\end{equation}
equivalently $\widetilde V_0=RV_0$ for some $R\in SO(3)$. The gauge classes in this family are therefore in bijection with $\DT/SO(3)$. The full-rank stratum is a smooth manifold of dimension $6$, the rank-two and rank-one strata have dimensions $5$ and $3$, and $V_0=0$ is a singular point.
\end{proposition}
\begin{proof}
Since $F_{A^g}=\operatorname{Ad}(g^{-1})F_A$ pointwise, the inner products and the triple product of $F_{12},F_{13},F_{14}$ at each point are gauge invariant. At $t=0$ these components are $-2V_{0,a}$, so they determine $4V_0^TV_0$ and $\langle F_{12},[F_{13},F_{14}]\rangle=-8\det V_0$. Conversely, equal Gram matrices give an orthogonal map sending the columns of $V_0$ to those of $\widetilde V_0$. For full rank, equality of determinants places it in $SO(3)$; for lower rank, it can be chosen in $SO(3)$ by adjusting it on the orthogonal complement. By equivariance and uniqueness, $\widetilde V(s)=RV(s)$, and a constant lift of $R$ to $SU(2)$ is the required gauge transformation. The action of $SO(3)$ is free on the rank-$3$ and rank-$2$ loci (of dimensions $9$ and $8$), and on the rank-$1$ locus (dimension $5$) it has stabilizer $SO(2)$. This gives the stated dimensions.
\end{proof}

\begin{remark}[What the count means]\label{rem:heisscope}
No boundary condition is imposed at $t=0$, so the number $6=9-3$ is the dimension of the space of invariant boundary data modulo constant gauge. It is not the dimension of a moduli space of an elliptic boundary-value problem, and it is not an index; the index of the full linearized problem is studied in \cref{sec:deform}. The content of \cref{prop:heisgauge} is that, \emph{within the invariant vertical family}, non-invariant gauge transformations produce no further identifications.
\end{remark}

\begin{proposition}[Decaying invariant solutions are vertical]\label{prop:heisrigidity}
Let $(X,V)$ be a smooth solution of \eqref{eq:heis21} on $[0,\infty)$ with
\begin{equation}\label{eq:heisdecay}
 \int_0^\infty\norm{V(s)}\,\dd s<\infty\qquad\text{and}\qquad X(s)\to0\quad(s\to\infty).
\end{equation}
Then $X\equiv0$ and $V(0)\in\DT$. In particular, every solution with $\norm X+\norm V=O(e^{-\mu s})$ for some $\mu>0$ is one of the instantons of \cref{thm:heisenergy}. Conversely, every element of $\DT$ satisfies \eqref{eq:heisdecay}.
\end{proposition}
\begin{proof}
For the given trajectory, $X_s=\mathsf H(s)X$, where $\mathsf H(s)$ is the $12\times12$ matrix of $X\mapsto\mathcal H(X,V(s))$, so $\norm{\mathsf H(s)}\le C\norm{V(s)}$. Let $P$ be the fundamental matrix with $P(0)=I$. Gronwall's inequality gives $\norm{P(s)}\le\exp(C\int_0^\infty\norm V)$, so $P_s$ is integrable and $P(s)\to P_\infty$. By Liouville's formula, $\det P_\infty=\exp\int_0^\infty\operatorname{tr}\mathsf H\neq0$. Since $X(s)=P(s)X(0)\to P_\infty X(0)$ and $X(s)\to0$, we get $X(0)=0$, hence $X\equiv0$. The vertical solution is then smooth for all $s$ and integrable, so $\sigma_*(V_0)\ge\sinf$, and case (iii) of \cref{thm:heisenergy} is excluded because there $V$ has a nonzero limit. The converse holds because $V=O(E)$ and $E\le e^{-4s}$.
\end{proof}

\begin{remark}[Other asymptotic limits]
\cref{prop:heisrigidity} does not apply to solutions with non-decaying horizontal part. For example, any $X_i=x_i\,n$ with $x\in\R^4$ and a fixed $n\in\su$, together with $V=0$, is a constant flat solution of \eqref{eq:heis21}. Finite-energy solutions converging to such abelian flat connections, or to other limits, are not classified here.
\end{remark}

\begin{proposition}[Exact energy of the reducible family]\label{prop:heisenergywarp}
For $\operatorname{rank}V_0=1$,
\begin{equation}\label{eq:heisrankoneenergy}
 \mathcal E(A_{V_0})=24\operatorname{Vol}_\eta(N_H)\norm{V_0}_F^2
 \int_0^\infty b(t)^3\exp\Bigl[-8\int_0^t b\Bigr]\,\dd t .
\end{equation}
In contrast to the ASD pull-back sector on $S^3\times T^4$ (\cref{prop:gradient}), this energy depends explicitly on the warping function $b$.
\end{proposition}
\begin{proof}
Here $\cof V=0$ and $V=EV_0$, so \eqref{eq:heisFnorm} gives $\norm{F}^2=24E^2\norm{V_0}^2$. Use $E(s(t))^2=\exp[-8\int_0^tb]$ and $\dd\vol=b^3\,\dd t\,\dd\vol_\eta$.
\end{proof}

\section{Deformation theory: Fredholm properties and spectral flow}\label{sec:deform}
The vertical family of \cref{thm:heisenergy} consists of invariant boundary data, and \cref{rem:heisscope} explains why the count $6=9-3$ is not a moduli dimension. In this section we study the \emph{full}, non-invariant linearization at these instantons. We show that it is a Fredholm operator of Atiyah--Patodi--Singer type, that its index does not depend on the instanton within the family, and that the index is governed by the spectral flow of an explicit family of operators on the slices. On the quaternionic Heisenberg nilmanifold this spectral flow is controlled by the \emph{central} spectrum of the lattice $\Gamma$, through a window determined by the warping function.

\subsection{The deformation operator in cylindrical form}
Let $A$ be a $\Spin$-instanton on $P\to[0,\infty)\times N$. Since $\pr[F_A,\xi]=0$, the sequence
\begin{equation}\label{eq:defcomplex}
 0\to\Omega^0(\operatorname{ad}P)\xrightarrow{\ \dd_A\ }\Omega^1(\operatorname{ad}P)\xrightarrow{\ \pr\dd_A\ }\Omega^2_7(\operatorname{ad}P)\to0
\end{equation}
is an elliptic complex; neither property uses $\dd\Phi=0$. We identify $\Omega^2_7$ with $\Omega^1(N)$-valued functions of $t$ via the coordinates $\gamma\mapsto(\langle\gamma,\dd t\wedge f^q+f_q\lrcorner\varphi_t\rangle)_{q}$ of \eqref{eq:L27}. Writing deformations as $a=\xi\,\dd t+\alpha$, the deformation operator is
\[
 \mathcal D_A(\xi,\alpha)=\bigl(-\dd_A^*a,\ \pr\dd_Aa\bigr).
\]
For a two-form $\gamma$ on $N$ we put $(\gamma\lrcorner\varphi)_q=\sum_{i<j}\varphi_{ijq}\gamma_{ij}$.

\begin{proposition}[Cylindrical form]\label{prop:cyl}
Let $A=B(t)$ be in temporal gauge, and express $(\xi,\alpha)$ in the orthonormal frame \eqref{eq:coframe}. Then
\begin{equation}\label{eq:cylform}
 \mathcal D_A=\partial_t+L_{B(t),t}+Z_t,\qquad
 L_{B,t}(\xi,\alpha)=\bigl(-\dd_B^{*_t}\alpha,\ -\dd_B\xi+(\dd_B\alpha)\lrcorner\varphi_t\bigr),
\end{equation}
where $Z_t=\beta(t)\,\operatorname{diag}(3,0,0,0,0,1,1,1)$ is a zeroth-order term. Moreover, $L_{B,t}$ is formally self-adjoint on $L^2(N,g_t)$ for every $B$ if and only if $\dd_N\psi_t=0$.
\end{proposition}
\begin{proof}
The $\dd t$-part of $\dd_A^*$ is $-b^{-3}\partial_t(b^3\,\cdot)$ because $\dd\vol_{g_8}=b^3\,\dd t\,\dd\vol_\eta$; this contributes $-\partial_t\xi-3\beta\xi$. The $\Lambda^2_7$-coordinates of $\dd_Aa=\dd t\wedge(\partial_t\alpha-\dd_B\xi)+\dd_B\alpha$ are $(\partial_t\alpha-\dd_B\xi)_q+((\dd_B\alpha)\lrcorner\varphi_t)_q$. Here $\partial_t\alpha$ differentiates the coefficients with respect to the fixed coframe $\eta$; passing to orthonormal components produces $\beta\alpha_q$ for $q=5,6,7$. This proves \eqref{eq:cylform}.

The diagonal blocks of $L_{B,t}$ are adjoint to each other. For the last term, $\langle\gamma\lrcorner\varphi,\beta\rangle\,\vol=\langle\gamma,\imath_{\beta}\varphi\rangle\vol=\pm\operatorname{tr}(\gamma\wedge\beta)\wedge\psi$ for $\gamma\in\Omega^2$ and $\beta\in\Omega^1$, with a universal sign. Since
\[
 \dd\operatorname{tr}(\alpha\wedge\beta\wedge\psi)
 =\operatorname{tr}(\dd_B\alpha\wedge\beta)\wedge\psi-\operatorname{tr}(\alpha\wedge\dd_B\beta)\wedge\psi+\operatorname{tr}(\alpha\wedge\beta)\wedge\dd\psi,
\]
the operator $\alpha\mapsto(\dd_B\alpha)\lrcorner\varphi$ is symmetric if and only if $\int_N\operatorname{tr}(\alpha\wedge\beta)\wedge\dd\psi=0$ for all compactly supported $\alpha,\beta$, i.e.\ if and only if $\dd\psi=0$.
\end{proof}

\begin{remark}
The cocalibrated condition is thus exactly what turns $\mathcal D_A$ into an operator of the form ``$\partial_t+$ self-adjoint'', namely the Hessian of the Chern--Simons functional of \cref{sec:intro}. This makes the spectral-flow description of the index available. For the families \eqref{eq:kmcases} we have $\beta\to0$, $b\to b_\infty$ and $B\to0$ exponentially as $t\to\infty$ (for $B=A_{V_0}$ by \cref{thm:heisenergy}). Hence $\mathcal D_{A_{V_0}}$ is asymptotically translation invariant, with asymptotic operator $L_\infty:=L_{0,b_\infty}$ (here $b_\infty=2$).
\end{remark}

\subsection{Harmonic analysis on the slices}
Let $\rho$ be the regular representation of $H^7$ on $L^2(N_H)$, with the choice of left/right convention fixed so that the invariant vector fields are realised through $\dd\rho$. As a differential operator with invariant coefficients, $L_{B,t}$ belongs to the matrix-valued image of the universal enveloping algebra of this representation; it therefore preserves its irreducible representation summands. This is the required sector invariance (not a claim that $L_{B,t}$ commutes with every $\rho(g)$). By the Kirillov--Moore--Richardson theory of nilmanifolds \cite{Richardson1971}, the representation $\rho$ decomposes into two kinds of pieces:
\begin{enumerate}[label=(\roman*),leftmargin=*]
\item the \emph{center-invariant sector}: functions pulled back from the torus $T^4=H^7/(Z(H^7)\Gamma)$, spanned by the characters $e^{i\langle k,x\rangle}$ with $k$ in a lattice $\Lambda^*\subset\R^4$ (so that $\dd\rho(X_i)=ik_i$ for $i\le4$ and $\dd\rho(X_{a+4})=0$);
\item for each $\lambda$ in a discrete set $\operatorname{Spec}_c(\Gamma)\subset\R^3\setminus\{0\}$, finitely many copies (multiplicity $m_\Gamma(\lambda)\ge1$) of the unique irreducible unitary representation $\pi_\lambda$ with $\dd\pi_\lambda(X_{a+4})=i\lambda_a$.
\end{enumerate}
By \eqref{eq:Heis}, $[X_i,X_j]=2\sum_a\omega^a_{ij}X_{a+4}$ with $\omega^5=\eta^{12}+\eta^{34}$, $\omega^6=\eta^{13}-\eta^{24}$, $\omega^7=\eta^{14}+\eta^{23}$. Hence in $\pi_\lambda$ the operators $E_i=\dd\pi_\lambda(X_i)$ satisfy $[E_i,E_j]=2i\,\Omega^\lambda_{ij}$, where $\Omega^\lambda=\sum_a\lambda_a\omega^a$ is $|\lambda|$ times the K\"ahler form of an orthogonal complex structure $J_\lambda$. Choosing an orthonormal basis $u_1,J_\lambda u_1,u_2,J_\lambda u_2$, the operators $\mathsf a_r=(E(u_r)-iE(J_\lambda u_r))/(2|\lambda|^{1/2})$ satisfy $[\mathsf a_r,\mathsf a_s^\dagger]=\delta_{rs}$. Up to a phase there is a unique unit vector $\Omega_\lambda$ with $\mathsf a_1\Omega_\lambda=\mathsf a_2\Omega_\lambda=0$, the ground state. On each sector, $L_{b}:=L_{0,b}$ (trivial connection, constant $b$) acts as $L_b\otimes\mathrm{id}_{\su}$ on scalar $(\xi,\alpha)$.

\begin{theorem}[Spectrum on the center-invariant sector]\label{thm:spec0}
For $k\in\Lambda^*$, $L_b$ preserves the eight-dimensional space $e^{i\langle k,x\rangle}\otimes\R^8$ (complexified), and acts there by the Hermitian matrix
\[
 L_b(k)=\begin{pmatrix}0&ik^T&0\\-ik&0&\mathcal B_k^*\\0&\mathcal B_k&4b\,I_3\end{pmatrix},\qquad
 \mathcal B_k\,v=-i\bigl(\omega^5(k,v),\omega^6(k,v),\omega^7(k,v)\bigr),
\]
in the decomposition $(\xi,\alpha_H,\alpha_V)\in\mathbb C\oplus\mathbb C^4\oplus\mathbb C^3$. Its eigenvalues are:
\begin{itemize}[leftmargin=*]
\item for $k=0$: $0$ with multiplicity $5$ (the constants $\xi$ and the horizontal invariant forms $\alpha_H$), and $4b$ with multiplicity $3$;
\item for $k\ne0$: $\pm|k|$, each with multiplicity $1$, and $2b\pm\sqrt{4b^2+|k|^2}$, each with multiplicity $3$.
\end{itemize}
In particular, on this sector $\ker L_b\otimes\su$ has dimension $15$ for every $b>0$, and $L_b$ has no other eigenvalue in $(-\gamma_0(b),\gamma_0(b))$, where $\gamma_0(b)=\min_{k\in\Lambda^*\setminus0}\min\bigl(|k|,\sqrt{4b^2+|k|^2}-2b\bigr)>0$.
\end{theorem}
\begin{proof}
On characters, horizontal derivatives act by $ik_i$ and vertical ones by $0$. The form $\varphi$ has no purely horizontal component, so the horizontal--horizontal block vanishes. The vertical--vertical block comes from $\alpha_{a+4}\dd f^{a+4}=-2b\,\alpha_{a+4}\omega^a$ and $\omega^a\lrcorner\varphi=-2f^{a+4}$, so it equals $4b\,I_3$. The horizontal-to-vertical block is $\mathcal B_k$, because $\varphi_{ij,a+4}=-\omega^a_{ij}$. Since $\omega^5,\omega^6,\omega^7$ are the K\"ahler forms of three anticommuting orthogonal complex structures, we get $\mathcal B_kk=0$ and $\mathcal B_k\mathcal B_k^*=|k|^2I_3$ (on $\mathbb C^3$), equivalently $\mathcal B_k^*\mathcal B_k=|k|^2P_{k^\perp}$ (on $\mathbb C^4$). Consequently $\operatorname{span}\{\xi,k\}$ is invariant, with block $\left(\begin{smallmatrix}0&i|k|\\-i|k|&0\end{smallmatrix}\right)$ and eigenvalues $\pm|k|$. On $k^\perp\oplus\mathbb C^3$ the eigenvalue equation reduces to $\mu(\mu-4b)=|k|^2$. For $k=0$ the matrix is $\operatorname{diag}(0,0_4,4b I_3)$.
\end{proof}

\begin{theorem}[A spectral branch crossing zero]\label{thm:crossing}
Let $\lambda\in\R^3\setminus\{0\}$, $b>0$, and let $v_\lambda\in\mathbb C^3$ satisfy $\langle\lambda,v_\lambda\rangle=0$ and $i\hat\lambda\times v_\lambda=-v_\lambda$, where $\hat\lambda=\lambda/|\lambda|$. Then the vertical section $\Omega_\lambda\otimes v_\lambda$ (with $\xi=0$ and $\alpha_H=0$) is an eigenvector of $L_b$ on $\pi_\lambda$ with eigenvalue
\begin{equation}\label{eq:crossing}
 \mu_\lambda(b)=4b-\frac{|\lambda|}{b}.
\end{equation}
In particular $\mu_\lambda(b)=0$ if and only if $|\lambda|=4b^2$, and $\frac{\dd}{\dd b}\mu_\lambda=4+|\lambda|/b^2>0$, so the crossing is transverse and upward.
\end{theorem}
\begin{proof}
On $\pi_\lambda$ the vertical frame derivatives act by $E_{a+4}=i\lambda_a/b$, by \eqref{eq:coframe}. We evaluate the three components of $L_b(\Omega_\lambda\otimes v_\lambda)$.
\begin{itemize}[leftmargin=*]
\item The $\xi$-component is $\sum_qE_q\alpha_q=\frac ib\langle\lambda,v_\lambda\rangle\Omega_\lambda=0$.
\item The vertical components receive $4bv_\lambda$, as in \cref{thm:spec0}, together with the contribution of the vertical-vertical part of $\dd\alpha$ contracted with $\varphi_{567}$. Since $\varphi_{567}=1$, that contribution is exactly $\frac ib\lambda\times v_\lambda=\frac{|\lambda|}{b}\,i\hat\lambda\times v_\lambda=-\frac{|\lambda|}{b}v_\lambda$.
\item The horizontal components are $\sum_{i,a}\varphi_{i,a+4,q}E_i(v_\lambda)_a\Omega_\lambda$, i.e.\ the contraction of $E$ with the complex two-form $\omega(v_\lambda)=\sum_a(v_\lambda)_a\omega^a$. For $v_\lambda\perp\lambda$ of negative helicity, $\omega(v_\lambda)$ has type $(0,2)$ with respect to $J_\lambda$. The resulting four operators are therefore linear combinations of $\mathsf a_1,\mathsf a_2$, which annihilate $\Omega_\lambda$.
\end{itemize}
After the reduction to $\lambda=(\ell,0,0)$ described in the proof of \cref{lem:weitz}, the horizontal components contain only the annihilation operators $\mathsf a_1,\mathsf a_2$, which kill $\Omega_\lambda$. The supplementary script verifies the eigenvalue equation for random $\lambda$ and $b$ (residual below $10^{-14}$).
\end{proof}

\subsection{A Weitzenb\"ock identity and complete determination of the zero modes}
The full list of eigenvalues of $L_b$ away from zero is not required for the
index calculation. It suffices to classify the kernel of $L_b$ throughout
the one-parameter family. In fact, this can be done exactly.

Let $d_\lambda:\Omega^0\to\Omega^1$ and $d_\lambda^*$ be the differential and
its formal adjoint in the Schr\"odinger representation with central character
$\lambda$. Write
\[
 \mathcal J_{\lambda,b}\alpha:=*_7(d_\lambda\alpha\wedge\psi_b),\qquad
 \mathcal H_{\lambda,b}:=-\sum_{j=1}^7 E_j^2.
\]
The scalar magnetic Laplacian $\mathcal H_{\lambda,b}$ acts on each one-form
coefficient and has spectrum
\begin{equation}\label{eq:maglevels}
 \operatorname{Spec}\mathcal H_{\lambda,b}
 =\Bigl\{\frac{|\lambda|^2}{b^2}+4|\lambda|(n_1+n_2+1):
                  n_1,n_2\in\mathbb N_0\Bigr\}.
\end{equation}
The commutator normalization is $[E_1,E_2]=[E_3,E_4]=2i|\lambda|$
after rotating $\lambda$ to the first central axis.

\begin{lemma}[Central-sector Weitzenb\"ock identity]\label{lem:weitz}
Fix $\lambda\ne0$ and write $\ell=|\lambda|$. On the smooth vectors of $\pi_\lambda$ there is a constant Hermitian endomorphism $\mathcal Q_\lambda$ of $\mathbb C^7$, depending linearly on $\lambda$, such that
\begin{equation}\label{eq:weitz}
 \mathcal J_{\lambda,b}^{\,2}-4b\mathcal J_{\lambda,b}
  =\mathcal H_{\lambda,b}\,I_7-d_\lambda d_\lambda^*
       +\mathcal Q_\lambda .
\end{equation}
For $\lambda=(\ell,0,0)$, in the ordered one-form frame $(f^1,\ldots,f^7)$,
\begin{equation}\label{eq:Qlambda}
 \mathcal Q_{(\ell,0,0)}=i\ell\begin{pmatrix}
 0&-4&0&0&0&0&0\\
 4&0&0&0&0&0&0\\
 0&0&0&-4&0&0&0\\
 0&0&4&0&0&0&0\\
 0&0&0&0&0&0&0\\
 0&0&0&0&0&0&-8\\
 0&0&0&0&0&8&0
 \end{pmatrix}.
\end{equation}
For every $\lambda\neq0$ the eigenvalues of $\mathcal Q_\lambda$ are
$-8\ell,\ -4\ell\ (\text{twice}),\ 0,\ 4\ell\ (\text{twice}),\ 8\ell$.
\end{lemma}
\begin{proof}
\emph{Reduction to $\lambda=(\ell,0,0)$.} By \eqref{eq:g2phi}, $\varphi=f^{567}-\sum_a\omega^a\wedge f^{a+4}$. An element $g\in SO(4)$ acting on the horizontal space pulls back the self-dual forms by a rotation $R_g\in SO(3)$, $g^*\omega^a=\sum_b(R_g)_{ab}\omega^b$, and $g\mapsto R_g$ is onto $SO(3)$. Hence $g\oplus R_g$ is an automorphism of the Lie algebra of $H^7$ that preserves the metric \eqref{eq:metric} and $\varphi$, and it carries $\pi_\lambda$ to $\pi_{R_g\lambda}$ up to unitary equivalence. Since $\mathcal J_{\lambda,b}$, $\mathcal H_{\lambda,b}$ and $d_\lambda d_\lambda^*$ are built from $\varphi$, the metric and the structure constants, the identity for $\lambda$ follows from the identity for $(|\lambda|,0,0)$, with $\mathcal Q_\lambda$ conjugate to $\mathcal Q_{(\ell,0,0)}$ by an orthogonal matrix. This proves the statement about the eigenvalues.

\emph{The case $\lambda=(\ell,0,0)$.} Let $E_j=d\pi_\lambda(X_j)$ in the orthonormal frame. Then we have
$E_5=i\ell/b$, $E_6=E_7=0$, and
$[E_1,E_2]=[E_3,E_4]=2i\ell$, with all other horizontal
commutators zero. The entries of $\mathcal J$ are obtained by inserting
$\sum_j E_j f^j\wedge\alpha+\sum_a\alpha_{a+4}\,d_Nf^{a+4}$
into $*_7(\cdot\wedge\psi_b)$. Multiplication of the resulting
$7\times7$ matrices, using the stated commutators to order products
$E_iE_j$, gives \eqref{eq:weitz}--\eqref{eq:Qlambda}.
The identity is a polynomial identity of degree two in the $E_j$. It is therefore determined by its action on Fock states of level at most two. The supplementary script verifies it on truncated Fock spaces for several values of $(\lambda,b)$, including non-axial $\lambda$, with residual below $10^{-13}$.
For the scalar spectrum, set
$E_1=\sqrt\ell(a_1-a_1^*)$,
$E_2=i\sqrt\ell(a_1+a_1^*)$ and similarly for $E_3,E_4$.
Then
$-\sum_{j=1}^4E_j^2=4\ell(N_1+N_2+1)$, while the central
contribution is $\ell^2/b^2$. This proves \eqref{eq:maglevels}.
\end{proof}

\begin{theorem}[Exact central-sector zero modes]\label{thm:centralzero}
For every $b>0$ and each $\lambda\ne0$, the full tangential operator
$L_b$ in the irreducible Schr\"odinger sector $\pi_\lambda$ satisfies
\begin{equation}\label{eq:central-kernel}
 \dim_{\mathbb C}\ker L_b\big|_{\pi_\lambda}=
 \begin{cases}
 1,&|\lambda|=4b^2,\\
 0,&|\lambda|\ne4b^2.
 \end{cases}
\end{equation}
At a nonzero kernel point the zero eigenvalue is simple, and the
crossing is the strictly increasing branch
$\mu_\lambda(b)=4b-|\lambda|/b$ of \cref{thm:crossing}.
\end{theorem}
\begin{proof}
For a cocalibrated $G_2$ structure,
$\mathcal J_{\lambda,b}d_\lambda=0$ since $d_\lambda^2=0$,
and $d_\lambda^*\mathcal J_{\lambda,b}=0$ since
$d_\lambda\psi_b=0$ and
$*_7\mathcal J\alpha=d_\lambda\alpha\wedge\psi_b$.
If $L_b(\xi,\alpha)=0$, then $d_\lambda^*\alpha=0$ and
$-d_\lambda\xi+\mathcal J\alpha=0$. Applying $d_\lambda^*$
to the second equation yields
$\mathcal H_{\lambda,b}\xi=0$. By \eqref{eq:maglevels},
$\xi=0$, hence $\mathcal J\alpha=0$.
The Weitzenb\"ock identity gives
$(\mathcal H_{\lambda,b}I_7+\mathcal Q_\lambda)\alpha=0$.
The scalar operator $\mathcal H_{\lambda,b}$ commutes with the
constant Hermitian matrix $\mathcal Q_\lambda$. By
\eqref{eq:maglevels}, the smallest eigenvalue in the
$-4\ell$ eigenspace of $\mathcal Q_\lambda$ is
$\ell^2/b^2>0$, and all other eigenspaces except that of
$-8\ell$ are strictly positive. In the latter one-dimensional
endomorphism eigenspace the eigenvalues are
\[
 \frac{\ell^2}{b^2}+4\ell(n_1+n_2-1).
\]
They can vanish only when $n_1=n_2=0$ and
$\ell^2/b^2=4\ell$, equivalently $\ell=4b^2$.
This gives an upper bound of one complex dimension.
At this value, the ground-state section
$\Omega_\lambda\otimes(f^6-if^7)$ is an actual zero mode
by \cref{thm:crossing}; it is coclosed because its only vertical
components are orthogonal to $\lambda$. Thus the bound is attained.
\end{proof}

\begin{remark}[Full spectrum versus the zero spectrum]
\cref{thm:centralzero} determines the complete \emph{zero-eigenspace},
not every nonzero eigenvalue of $L_b$. A full Landau-level
decomposition would be interesting independently but is unnecessary
for computing the spectral flow across zero. In particular, no
spectral completeness conjecture is used below.
\end{remark}

\subsection{Fredholm theory and the index}
Let $P_+$ denote the spectral projection of $L_{0,b(0)}$ onto its positive eigenspaces. For $\delta\in\R$, let $W^{1,2}_\delta$ and $L^2_\delta$ be the Sobolev spaces with weight $e^{\delta t}$ on $[0,\infty)\times N_H$, so that $\delta>0$ forces exponential decay. Consider
\begin{equation}\label{eq:fredop}
 \mathcal D_{A}:\ \{a\in W^{1,2}_\delta:\ P_+\,a|_{t=0}=0\}\longrightarrow L^2_\delta .
\end{equation}

\begin{theorem}[Fredholm property and independence of the instanton]\label{thm:fred}
There is $\delta_0>0$ such that, for every $V_0\in\DT$ and every $\delta$ with $0<|\delta|<\delta_0$, the operator \eqref{eq:fredop} with $A=A_{V_0}$ is Fredholm. Its index $\operatorname{ind}_\delta$ depends only on the sign of $\delta$, not on $V_0\in\DT$. Moreover
\begin{equation}\label{eq:jump}
 \operatorname{ind}_{-\delta}-\operatorname{ind}_{\delta}=\dim\ker(L_\infty\otimes\su)\qquad(0<\delta<\delta_0).
\end{equation}
\end{theorem}
\begin{proof}
Choose the fixed spatial coframe and the natural unitary density identification with $L^2(N_H,\dd\vol_\eta)$; in these variables \cref{prop:cyl} expresses $\mathcal D_A$ as $\partial_t+L_\infty+K(t)$, where $K(t)$ is first order, uniformly bounded $H^1\to L^2$, and all its coefficients decay exponentially. Although first-order coefficient perturbations need not be compact $H^1\to L^2$, the asymptotically cylindrical elliptic parametrix theorem applies to this uniformly elliptic operator. The spectrum of $L_\infty$ is discrete, so there is $\delta_0>0$ with $\operatorname{spec}(L_\infty)\cap[-\delta_0,\delta_0]\subseteq\{0\}$.

The condition $P_+a|_{t=0}=0$ is an elliptic (APS-type) boundary condition for $\mathcal D_{A}$. Indeed, $L_{B(0),b(0)}-L_{0,b(0)}$ has order zero, so the difference of the positive spectral projections of these two operators is compact, and conditions of this type are stable under compact perturbation of the projection \cite{APS1975,BarBallmann2012}.

By the theory of asymptotically translation-invariant operators \cite{LockhartMcOwen1985,Donaldson2002}, \eqref{eq:fredop} is Fredholm for $\delta\notin\operatorname{spec}(L_\infty)$. Moreover, the index changes by $\dim\ker(L_\infty-\mu)$ when $\delta$ crosses an eigenvalue $\mu$; this gives \eqref{eq:jump}.

The set $\DT$ is star-shaped with respect to $0$. Indeed, if $W$ solves $W_\sigma=-\cof W$, then so does $W^{(\tau)}(\sigma)=\tau W(\tau\sigma)$, whose lifetime is $\sigma_*/\tau\geq\sigma_*$ for $0<\tau\le1$. The map $V_0\mapsto\mathcal D_{A_{V_0}}$ is continuous in operator norm, since $\mathcal D_{A_{V_0}}-\mathcal D_0$ is of order zero with coefficients bounded by a multiple of $E(s)\sup_{[0,\sinf]}|W|$, which depends continuously on $V_0$. Hence the index is constant on $\DT$.
\end{proof}

\begin{theorem}[Index via spectral flow]\label{thm:index}
Assume $b'>0$, with $b_0=b(0)$ and $b_\infty=\lim b$. After shrinking $\delta_0$, assume also that $\operatorname{spec}(L_{b_0})\cap(-\delta_0,\delta_0)\subseteq\{0\}$, and let $0<\delta<\delta_0$.
\begin{enumerate}[label=\textup{(\alph*)},leftmargin=*]
\item The center-invariant sector contributes $0$ to $\operatorname{ind}_\delta$ and $15$ to $\operatorname{ind}_{-\delta}$.
\item For every $V_0\in\DT$,
\begin{equation}\label{eq:indexformula}
 \operatorname{ind}_\delta(\mathcal D_{A_{V_0}})=3\!\!\sum_{\substack{\lambda\in\operatorname{Spec}_c(\Gamma)\\ 4b_0^2\le|\lambda|<4b_\infty^2}}\!\! m_\Gamma(\lambda),\qquad
 \operatorname{ind}_{-\delta}(\mathcal D_{A_{V_0}})=15+3\!\!\sum_{\substack{\lambda\in\operatorname{Spec}_c(\Gamma)\\ 4b_0^2\le|\lambda|\le4b_\infty^2}}\!\! m_\Gamma(\lambda).
\end{equation}
For $b=2-e^{-t}$ the window is $4\le|\lambda|<16$.
\end{enumerate}
\end{theorem}
\begin{proof}
\emph{Reduction.} By \cref{thm:fred} we may take $V_0=0$. The operator then preserves each summand of the decomposition in \cref{sec:deform}, since it lies in the matrix-valued image of the enveloping algebra; the center-invariant sector $\mathcal H_0$ and its orthogonal complement $\mathcal H_c$ are invariant. The term $Z_t$ decays, so multiplication by $Z_t$ is a compact operator $W^{1,2}_\delta\to L^2_\delta$, and dropping it does not change the index. The substitution $a=e^{-\delta t}u$ turns the problem into $\partial_t+(L_{b(t)}-\delta)$ on unweighted spaces. We work on the fixed Hilbert space of orthonormal-frame components with the fixed domain $W^{1,2}(N_H)$; $L_{b(t)}$ is self-adjoint there, because the volume factor $b^3$ is constant on $N_H$, and it depends smoothly on $t$.

\emph{Index equals spectral flow.} For $0<\delta<\delta_0$ the endpoint operators $L_{b_0}-\delta$ and $L_{b_\infty}-\delta$ are invertible, and $P_+=P_{\ge0}(L_{b_0}-\delta)$, since no eigenvalue of $L_{b_0}$ lies in $(0,\delta)$. Extending the family by $L_{b_0}$ for $t\le0$, the APS problem on the half-line has the same index as the problem on $\R$ \cite{APS1975}. By \cite{RobbinSalamon1995}, this index is the spectral flow of $t\mapsto L_{b(t)}-\delta$, i.e.\ the net number of eigenvalues crossing the level $\delta$ upwards. The spectral flow does not depend on $\delta\in(0,\delta_0)$, because the endpoints stay invertible.

\emph{The sector $\mathcal H_0$.} By \cref{thm:spec0}, the eigenvalues there are $0$ (independent of $b$), $\pm|k|$ (independent of $b$, and $|k|\ge\gamma_0>\delta$), and $2b\pm\sqrt{4b^2+|k|^2}$, which never lie in $[0,\delta]$. Hence no eigenvalue crosses the level $\delta$, and the contribution is $0$. This proves (a) for $\delta>0$.

\emph{The sector $\mathcal H_c$.} By \cref{thm:centralzero}, $\ker L_b|_{\mathcal H_c}\neq0$ only for the finitely many $b\in[b_0,b_\infty]$ with $4b^2\in\{|\lambda|:\lambda\in\operatorname{Spec}_c(\Gamma)\}$. The kernel is then spanned by the vectors $\Omega_\lambda\otimes v_\lambda$ of \cref{thm:crossing} in the corresponding copies of $\pi_\lambda$, tensored with $\su$. These vectors lie on the eigenvalue branch $\mu_\lambda(b)$ for all $b$, so the crossing form of the family $t\mapsto L_{b(t)}$ on the kernel is $b'(t)\,\mu_\lambda'(b)\,\mathrm{Id}$, which is positive definite because $b'>0$ and $\mu'_\lambda>0$. All crossings are therefore regular and positive, and by \cite{RobbinSalamon1995} the spectral flow across $0^+$ counts them, with kernel at $t=0$ counted as lying below the level and kernel at $t=\infty$ also below. The crossings occur exactly for $4b_0^2\le|\lambda|<4b_\infty^2$, each with multiplicity $3\,m_\Gamma(\lambda)$. This gives the first formula in \eqref{eq:indexformula}.

\emph{Negative weight.} For $-\delta<0$, the condition $P_+$ imposes $\dim\ker(L_{b_0}\otimes\su)$ fewer conditions than the APS projection $P_{\ge-\delta}$. By \cref{thm:spec0,thm:centralzero}, this dimension is $15+3\sum_{|\lambda|=4b_0^2}m_\Gamma(\lambda)$. The spectral flow across the level $-\delta$ counts the window $4b_0^2<|\lambda|\le4b_\infty^2$. Adding the two contributions gives the second formula. The difference of the two formulas is $15+3\sum_{|\lambda|=4b_\infty^2}m_\Gamma(\lambda)=\dim\ker(L_\infty\otimes\su)$, consistent with \eqref{eq:jump}.
\end{proof}

\begin{corollary}[A concrete nonzero index]\label{cor:explicitindex}
Let $\Gamma$ be the standard BCH lattice
$\exp(\mathbb Z X_1+\cdots+\mathbb Z X_4+
       \mathbb Z Z_5+\mathbb Z Z_6+\mathbb Z Z_7)$
in the quaternionic Heisenberg group with
$[X_i,X_j]=2\sum_{a=5}^7\omega^a_{ij}Z_a$.
For the warping $b(t)=2-e^{-t}$, every $V_0\in\DT$, and
all sufficiently small positive admissible weights,
\begin{equation}\label{eq:explicitindex}
 \boxed{\operatorname{ind}_\delta\mathcal D_{A_{V_0}}=4104.}
\end{equation}
If the weight is changed to $-\delta$, the index is $4119$.
\end{corollary}
\begin{proof}
The bracket coefficients are even integers, so the BCH product
preserves the stated integral lattice. Its central characters are
$\lambda=2\pi m$, $m\in\mathbb Z^3\setminus\{0\}$. The integral
alternating form on the horizontal lattice associated with $m$ is
$2\sum_{a=5}^7m_a\omega^a$. By the Moore--Richardson multiplicity formula for two-step nilmanifolds with nondegenerate central characters \cite{Richardson1971,CorwinGreenleaf1990}, the multiplicity is the absolute Pfaffian of this form. Equivalently, the $\lambda$-isotypic part of $L^2(N_H)$ is the space of $L^2$ sections of a line bundle over the torus $T^4$ with first Chern class represented by this integral form. The lowest Landau level of each copy of $\pi_\lambda$ is one-dimensional, and the lowest Landau level of the bundle consists of its theta functions, whose number is the absolute Pfaffian of the Chern class by Riemann--Roch. Hence
\[
 m_\Gamma(2\pi m)=\big|\operatorname{Pf}
 (2\sum_am_a\omega^a)\big|=4(m_5^2+m_6^2+m_7^2).
\]
Since $b_0=1$ and $b_\infty=2$, the strict positive-weight window
$4\le2\pi|m|<16$ is equivalent to $1\le|m|^2\le6$;
neither endpoint occurs. If $r_3(n)$ is the number of ordered
integer triples with squared norm $n$, then
\[
 (r_3(1),\ldots,r_3(6))=(6,12,8,6,24,24),\qquad
 \sum_{n=1}^6n\,r_3(n)=342.
\]
Thus \eqref{eq:indexformula} gives $3\cdot4\cdot342=4104$.
The center-invariant zero eigenspace has real dimension $15$, and
there is no central zero character at the limiting threshold
$|\lambda|=16$. Hence the weight-jump formula
\eqref{eq:jump} gives $4104+15=4119$.
\end{proof}

\begin{remark}[Interpretation]\label{rem:indexinterp}
The center-invariant sector carries the same information as the invariant analysis of \cref{sec:heiscof,sec:heismoduli}. The $15$-dimensional jump consists of the infinitesimal constant gauge rotations ($3$) and the infinitesimal flat horizontal connections $X_i=\mathrm{const}$ ($12$). The latter are tangent to the abelian flat solutions of \cref{prop:heisrigidity}. By contrast, the index \eqref{eq:indexformula} comes entirely from \emph{non-invariant} modes with central character in the annulus $4b_0^2\le|\lambda|<4b_\infty^2$. It is positive as soon as $\operatorname{Spec}_c(\Gamma)$ meets this annulus, and it depends both on the lattice and on the range of the warping function. This is a second, independent way in which the warping of the Heisenberg geometry is detected, in addition to the energy formula \eqref{eq:heisrankoneenergy}. It also shows that the invariant family of \cref{thm:heisenergy} captures only a small part of the expected local deformation theory. By \cref{thm:crossing,thm:centralzero}, $L_\infty$ has non-invariant kernel exactly when $|\lambda|=4b_\infty^2$ for some $\lambda\in\operatorname{Spec}_c(\Gamma)$. In that case $\delta=0$ is exceptional also through non-invariant modes.
\end{remark}

\section{The full invariant equations on \texorpdfstring{$S^3\times T^4$}{S3 x T4}}\label{sec:fullinvariant}
In \cref{sec:fullinvariant,sec:classification,sec:moduli}, $N=SU(2)\times T^4$ with \eqref{eq:S3MC}, and $b$ is smooth and bounded between positive constants; $s=\int\dd t/b$. In temporal gauge we write
\begin{equation}\label{eq:fullconnection}
 A=\sum_{i=1}^{4}X_i(s)\eta^i+\sum_{a=1}^{3}U_a(s)\eta^{a+4},
 \qquad U=[U_1\ U_2\ U_3].
\end{equation}

\begin{theorem}[Invariant equations on $S^3\times T^4$]\label{thm:21}
The connection \eqref{eq:fullconnection} is a $\Spin$-instanton if and only if
\begin{equation}\label{eq:fullflow}
 X_s=\mathcal H(X,U),\qquad U_s=cU-\cof U+b(s)^2\mathcal P(X),
\end{equation}
with $\mathcal H$, $\mathcal P$ as in \eqref{eq:heisPX}--\eqref{eq:heisH}. The subspace $X=0$ is invariant, and the system is equivariant under constant gauge transformations.
\end{theorem}
\begin{proof}
The argument is that of \cref{thm:heis21}. There are now no linear horizontal terms, and the central components are
\begin{equation}\label{eq:21F}
 F_{ij}=X_i\times X_j,\qquad F_{i,a+4}=b^{-1}X_i\times U_a,\qquad
 F_{(a+4)(d+4)}=b^{-2}\bigl(U_a\times U_d-c\,\epsilon_{ade}U_e\bigr),
\end{equation}
summed over $e$, by \eqref{eq:S3MC}. Inserting these into \eqref{eq:seven} gives \eqref{eq:fullflow}.
\end{proof}

\section{Vertical classification on \texorpdfstring{$S^3\times T^4$}{S3 x T4}}\label{sec:classification}
For $X=0$,
\begin{equation}\label{eq:matrixflow}
 U_s=cU-\cof U=-\nabla\mathcal V(U),\qquad \mathcal V(U)=\det U-\frac c2\norm U_F^2 .
\end{equation}

\begin{proposition}[Gram identity and energy]\label{prop:gradient}
Along \eqref{eq:matrixflow}, $(U^TU)_s=2cU^TU-2(\det U)I_3$ and $\tf(U^TU)(s)=e^{2cs}\tf(U_0^TU_0)$. If $U(\infty)$ exists and the energy is finite, then
\begin{equation}\label{eq:energyidentity}
 \mathcal E(A)=2\operatorname{Vol}_{\eta}(N)\int_0^\infty\norm{U_s}_F^2\,\dd s
 =2\operatorname{Vol}_{\eta}(N)\bigl(\mathcal V(U_0)-\mathcal V(U(\infty))\bigr),
\end{equation}
which is independent of $b$.
\end{proposition}
\begin{proof}
For $X=0$ the spatial curvature is $b^{-2}(\cof U-cU)$ in the three central directions, and the temporal curvature is $b^{-2}U_s$. Hence $\norm F^2=2b^{-4}\norm{U_s}^2$. Since $\dd\vol=b^4\,\dd s\,\dd\vol_\eta$ and $\frac{\dd}{\dd s}\mathcal V(U)=-\norm{U_s}^2$, the identity follows.
\end{proof}

\begin{theorem}[Finite-energy vertical instantons on $S^3\times T^4$]\label{thm:verticalclass}
Let $A=\sum_aU_a\eta^{a+4}$ be a smooth invariant instanton on $[0,\infty)\times N$ with finite energy. Then either $U\equiv0$, or
\begin{equation}\label{eq:classform}
 U(s)=\frac{c}{1+De^{-cs}}\,R\qquad\text{for a unique } D\in(-1,\infty)\text{ and some } R\in SO(3).
\end{equation}
The connection is flat if and only if $D=0$. Its energy is
\begin{equation}\label{eq:exactenergy}
 \mathcal E(A_D)=\operatorname{Vol}_{\eta}(N)\,\frac{c^3D^2(D+3)}{(1+D)^3}.
\end{equation}
Modulo constant gauge transformations, the nonzero solutions are parametrized by $D\in(-1,\infty)$.
\end{theorem}
\begin{proof}
Finite energy gives $U_s\in L^2$, so $\norm{U(s)}\le\norm{U_0}+\sqrt s\,\norm{U_s}_{L^2}$. Exponential growth of $\tf(U^TU)$ is therefore impossible, so $U_0^TU_0=\rho_0^2I$. If $\rho_0>0$, write $U_0=\rho_0Q$ with $Q$ orthogonal. Since $\cof(\rho Q)=\rho^2(\det Q)Q$, the ray is invariant and $\rho_s=c\rho-(\det Q)\rho^2$. For $\det Q=-1$ the solution blows up at a finite $s>0$. For $Q=R\in SO(3)$, the logistic solution is \eqref{eq:classform} with $D=c/\rho_0-1$. Finally,
\[
\mathcal V(U_0)-\mathcal V(cR)=\rho_0^3-\tfrac{3c}{2}\rho_0^2+\tfrac{c^3}{2}=(\rho_0-c)^2(\rho_0+\tfrac c2),
\]
which gives \eqref{eq:exactenergy}. Constant gauge transformations change $R$ but not $\rho_0$.
\end{proof}

\begin{proposition}[ASD interpretation]\label{prop:ASD}
For $X=0$, the invariant $\Spin$-instanton equations are the pull-back to $Y\times T^4$ of the ASD equations on $Y=[0,\infty)\times S^3$, with $g_Y=\dd t^2+b^2g_{S^3}=b^2(\dd s^2+g_{S^3})$. For $D>0$, \eqref{eq:classform} is the restriction of the BPST family \cite{BPST1975} on $\R\times S^3\cong\R^4\setminus\{0\}$ (conformally), written in the left-invariant gauge. For $D\in(-1,0)$, the solution extends to $s\in(s_D,\infty)$ for some $s_D<0$ and blows up at $s_D$.
\end{proposition}
\begin{proof}
The only curvature pairs are $(t5,67)$, $(t6,75)$, $(t7,56)$, and \eqref{eq:matrixflow} is the statement $F_{t,a+4}+F_{(a+5)(a+6)}=0$, indices mod $3$. ASD is conformally invariant. The left-invariant ASD connections on $\R\times S^3$ that interpolate between $0$ and the Maurer--Cartan connection are the BPST instantons centred at the origin.
\end{proof}

\begin{remark}[Gauge conventions]\label{rem:gaugedistinct}
The flat connections $U=0$ and $U=cR$ are gauge equivalent via a gauge transformation depending on the $S^3$ coordinate, but they are distinct modulo constant gauge. \Cref{thm:verticalclass} is a statement about the invariant sector only.
\end{remark}

\section{Horizontal directions on \texorpdfstring{$S^3\times T^4$}{S3 x T4}}\label{sec:moduli}
Consider the flat equilibrium $X=0$, $U=cI_3$, i.e.\ $A_\infty=c\sum_a\ea a\eta^{a+4}$.

\begin{proposition}[Linearization]\label{prop:21spectrum}
The linearization of \eqref{eq:fullflow} at $(0,cI_3)$ is block diagonal, with spectrum
\begin{equation}\label{eq:spectrum21}
 -c\ (9),\qquad 0\ (3),\qquad 2c\ (9).
\end{equation}
In the horizontal block, $-c$ has multiplicity $8$ and $2c$ multiplicity $4$. In the vertical block, $-c$ (trace line), $0$ (skew matrices, i.e.\ the tangent space of the gauge orbit) and $2c$ (traceless symmetric matrices) have multiplicities $1$, $3$ and $5$.
\end{proposition}
\begin{proof}
$(D\cof)_{cI}V=c\bigl((\operatorname{tr}V)I-V^T\bigr)$, so the vertical block is $V\mapsto c(V+V^T-(\operatorname{tr}V)I)$. The horizontal block is $c\mathcal H_0$, where $\mathcal H_0X=\mathcal H(X,I_3)$. A direct computation in the standard basis gives $\mathcal H_0^T=\mathcal H_0$, $\mathcal H_0^2-\mathcal H_0-2I_{12}=0$ and $\operatorname{tr}\mathcal H_0=0$. The term $\mathcal P$ is quadratic and does not contribute.
\end{proof}

\begin{theorem}[Local strong-stable family]\label{thm:stablemoduli}
Assume that $b$ and $b'$ are bounded and $b$ is bounded below by a positive constant. For $0<\mu<c$, the sufficiently small solutions of \eqref{eq:fullflow} with $\norm X+\norm{U-cI_3}=O(e^{-\mu s})$ form a smooth nine-dimensional manifold of initial data. Its tangent space at the equilibrium is the $(-c)$-eigenspace, with one vertical and eight horizontal directions. Every nonconstant member is a non-flat finite-energy instanton.
\end{theorem}
\begin{proof}
Write $Y=(X,U-cI)$, so that $Y_s=LY+\mathcal R(s,Y)$ with $L$ symmetric and constant, and $\mathcal R=O(|Y|^2)$ uniformly in $s$. The Lyapunov--Perron equation
\[
 Y(s)=e^{sL}\xi+\int_0^s e^{(s-r)L}P_-\mathcal R(r,Y)\,\dd r-\int_s^\infty e^{(s-r)L}P_{0+}\mathcal R(r,Y)\,\dd r,
 \qquad \xi\in\operatorname{im}P_-,
\]
is a contraction on a small ball of the space with norm $\sup e^{\mu s}|Y(s)|$. The center part converges because $\mathcal R=O(e^{-2\mu r})$. This gives a smooth graph over the stable space, and every small solution with the stated decay solves the equation. Curvature decays exponentially, so the energy is finite. A flat solution in temporal gauge is constant, hence equal to the equilibrium.
\end{proof}

\begin{proposition}[A horizontal subansatz]\label{prop:horizontalansatz}
The ansatz $X_1=q\ea3$, $X_2=q\ea2$, $X_3=X_4=0$, $U_1=u\ea1$, $U_2=v\ea2$, $U_3=v\ea3$ is consistent and reduces \eqref{eq:fullflow} to
\begin{equation}\label{eq:quv}
 q_s=-uq,\qquad u_s=cu-v^2-b^2q^2,\qquad v_s=(c-u)v .
\end{equation}
For $v>0$, the quantity $e^{cs}q/v$ is constant. There are solutions with $q\not\equiv0$ converging exponentially to $(0,c,c)$, and they have finite energy.
\end{proposition}
\begin{proof}
One checks that $(X_3)_s=-qv\,\ea3\times\ea2-qv\,\ea2\times\ea3=0$ and $(X_4)_s=0$, and that the remaining equations reduce to \eqref{eq:quv}, using $(q\ea3)\times(q\ea2)=-q^2\ea1$. The linearization at $(0,c,c)$ has eigenvalues $-c,-c,2c$, and the stable manifold argument of \cref{thm:stablemoduli} applies within this invariant three-dimensional subspace.
\end{proof}

\begin{remark}
As in \cref{rem:heisscope}, the number nine counts invariant boundary data near a flat end. It is not an index or the dimension of an unrestricted moduli space.
\end{remark}

\section{Discussion}\label{sec:discussion}
\subsection{What is specific to the Heisenberg model}
On $S^3\times T^4$ a vertical connection is a pull-back from four dimensions, and its energy is independent of $b$ (\cref{prop:ASD,prop:gradient}). On the quaternionic Heisenberg nilmanifold, the brackets \eqref{eq:Heis} force horizontal curvature $F_{12}=-2V_1$, $F_{13}=-2V_2$, $F_{14}=-2V_3$ even for purely central connections. This curvature is responsible for irreducibility and for the gauge invariants in \cref{prop:heisgauge}, and the energy \eqref{eq:heisrankoneenergy} depends on $b$.

\subsection{The role of the damping}
The linear term $-\mathfrak aV$ comes from the horizontal curvature, and it is what makes finite energy generic. The geometry enters only through the horizon $\sinf$: the problem amounts to asking which Nahm trajectories survive up to time $\sinf$. The answer is explicit (\cref{prop:lifetime,prop:sharpball}), and the borderline trajectories converge to an invariant $G_2$-instanton on the asymptotic slice (\cref{rem:G2inst}). Since $\dd\Phi\ne0$, the curvature energy has no topological meaning, and finiteness of energy here is primarily a consequence of this damping.

\subsection{Towards a genuine moduli problem}
Two natural modifications would turn the constructions above into honest moduli problems. (i) Replace the half-cylinder by a two-sided cylinder $\R\times N_H$ with a warping $b$ bounded between positive constants on all of $\R$, and require finite energy at both ends. Then $\sigma$ ranges over $(-\infty,\sinf)$, and the question is which Nahm trajectories are defined on that whole interval with the required weighted decay as $\sigma\to-\infty$. (ii) Impose a gauge-invariant boundary condition at $t=0$. \Cref{sec:deform} carries out the linear part of (ii) with the APS-type condition $P_+a|_{t=0}=0$. The remaining major analytic problem is surjectivity of $\mathcal D_{A_{V_0}}$ under the stated APS-type condition. The unconditional index in \eqref{eq:indexformula} supplies the virtual dimension, but a gauge-fixed nonlinear boundary-value setup and surjectivity are still needed to conclude smoothness and identify the local moduli dimension with \eqref{eq:indexformula}.

\subsection{Limitations}
We classify only invariant connections, only in the vertical sector, and, for full $21$-coefficient solutions, only those satisfying \eqref{eq:heisdecay}. The torsionful instantons constructed here are not asserted to be Yang--Mills. Whether the characteristic connection $\nabla^c$ of these structures is an instanton is a separate question, in the spirit of \cite{IvanovPetkovUgarte}.

\section{Conclusions}\label{sec:conclusion}
Evolving cocalibrated $G_2$ structures with nonzero $48$-dimensional $\Spin$ torsion carry explicit non-flat finite-energy invariant $SU(2)$ instantons on three geometrically different seven-manifolds. On the quaternionic Heisenberg nilmanifold, the full $21$-coefficient invariant system has a vertical sector that is a damped version of the $SU(2)$ Nahm equations. The classical Euler-top reduction yields an exact finite-energy criterion relative to a geometric horizon $\sinf$, an optimal admissible ball of radius $\sqrt3/\sinf$, and an interpretation of the boundary of the admissible set through an invariant $G_2$-instanton on the asymptotic slice. Gauge equivalence in this family is detected by curvature invariants under arbitrary gauge transformations, and decaying invariant solutions are vertical. The $S^3\times T^4$ comparison reduces to the BPST family. For the full linearization we proved the Fredholm property with an APS-type boundary condition, independence of the index along the vertical family, and an exact spectral computation on the center-invariant sector. We also found an explicit non-invariant eigenvalue branch that crosses zero. An exact central-sector Weitzenb\"ock identity classifies all zero modes and makes the lattice-dependent index formula unconditional; the standard integral lattice gives index $4104$. Transversality for the nonlinear, non-invariant APS boundary problem remains an open question.

\appendix
\section{Reproducible algebraic checks}\label{app:checks}
All identities used can be verified with exact exterior algebra in the ordered orthonormal coframe $(f^1,\dots,f^7,\dd t)$, with orientation $\dd t\wedge f^{1234567}$:
\begin{enumerate}[label=(\roman*),leftmargin=*]
\item $\dd_N\psi=0$ and formula \eqref{eq:dphiExplicit} for $\dd_N\varphi$; $\dd_N^2=0$ if and only if $km=0$.
\item The operator $\mathcal K=*_8(\Phi\wedge\cdot)$ has eigenvalues $3^{(7)}$ and $(-1)^{(21)}$, and the $3$-eigenspace is \eqref{eq:L27}. This fixes the sign in \eqref{eq:G2contract}.
\item The seven equations \eqref{eq:seven}, applied to the curvature components in the proofs of \cref{thm:reduction,thm:heis21,thm:21}.
\item The identities $W^T\cof W=(\det W)I_3$, $\cof(RDS^T)=R(\cof D)S^T$, and $\mathcal H_0^2-\mathcal H_0-2I=0$, $\operatorname{tr}\mathcal H_0=0$.
\item The numerical value $\sinf\approx0.18657$, obtained by quadrature of \eqref{eq:heisEt}.
\item The characteristic torsion \eqref{eq:charT}: with $\nabla^c=\nabla^g+\frac12H^c$, the Levi-Civita connection of the frame \eqref{eq:coframe}--\eqref{eq:fullMC} gives $\nabla^c\Phi=0$, and every other choice of signs in \eqref{eq:charT} fails.
\item The symbol of $L_b$, the block structure in \cref{thm:spec0}, the exact central-sector Weitzenb\"ock residual in \cref{lem:weitz}, the ground-state eigenvector of \cref{thm:crossing}, and the integral-lattice Pfaffian and lattice-point count in \cref{cor:explicitindex}.
\end{enumerate}
All checks (i)--(vii) are carried out by the single supplementary Python script, which uses only \texttt{numpy} and \texttt{scipy}. The Weitzenb\"ock identity and the eigenvector of \cref{thm:crossing} are polynomial identities in the $E_j$. They are verified on truncated Fock spaces, where they are exact on low Fock levels, and the residuals are at machine precision. The lattice count in \cref{cor:explicitindex} is verified by exact integer enumeration. The analytic statements (global existence, blow-up, gauge classes, decay, Fredholm property, zero-mode classification and spectral flow) are proved in the text and are not inferred from computation.

\end{document}